\documentclass[11pt]{article}
\usepackage[utf8]{inputenc}
\usepackage[margin=0.9in]{geometry}
\usepackage{amsmath,amssymb,amsthm}
\usepackage{graphicx}
\usepackage{booktabs}
\usepackage{multirow}
\usepackage{natbib}
\usepackage{xcolor}
\usepackage{caption}
\usepackage{subcaption}
\usepackage{hyperref}
\hypersetup{
    colorlinks,
    linkcolor={red!50!black},
    citecolor={blue!50!black},
    urlcolor={blue!80!black}
}
\usepackage{authblk}
\usepackage{tikz}
\usepackage{pgfplots}
\usepackage{verbatim}
\pgfplotsset{compat=1.18}
\usetikzlibrary{patterns,decorations.pathreplacing,calc,arrows.meta}
\usepackage{euler,charter}
\DeclareMathOperator{\Var}{Var}
\DeclareMathOperator{\diam}{diam}

\newcommand{\R}{\mathbb{R}}

\newcommand{\Alpha}{\mathrm{Alpha}}

\newtheorem{theorem}{Theorem}[section]

\newtheorem{proposition}[theorem]{Proposition}
\newtheorem{corollary}[theorem]{Corollary}

\newtheorem{remark}[theorem]{Remark}

\title{Detecting Regime Transitions in Dynamical Systems\\ via the Mixup Euler Characteristic Profile}

\author[1]{Sushovan Majhi\thanks{Corresponding Author: s.majhi@gwu.edu}}
\author[2]{Atish Mitra\thanks{amitra@mtech.ed}}
\author[3]{Santanu Nandi\thanks{santanu282@gmail.com}}
\author[4]{Md Nurujjaman\thanks{md.nurujjaman@nitsikkim.ac.in}}
\author[4]{Buddha Nath Sharma\thanks{bnsharma09@yahoo.com}}

\affil[1]{Data Science Program, George Washington University, Washington D.C., USA}
\affil[2]{Department of Mathematical Sciences, Montana Technological University, Butte, MT, USA}
\affil[3]{VIT-AP University, Amravati, India}
\affil[4]{Department of Physics, National Institute of Technology Sikkim, Ravangla, India}

\date{}

\begin{document}

\maketitle

\begin{abstract}
We develop a framework for detecting regime transitions in dynamical systems using the \emph{Mixup Euler Characteristic Profile} (Mixup ECP)---the Euler characteristic of the geometric intersection of ball unions around adjacent delay-embedded trajectory segments, viewed as a function of filtration scale. The Mixup ECP provides a detection statistic with a built-in null and guaranteed stability. We formalize regime detection as a low-side-permutation test, establish its validity and consistency, and introduce a multi-delay extension that automatically selects the most informative dynamical timescale. Complementing the topological signal with Complexity Variance, Higuchi fractal dimension, and a rolling mean baseline, the four-signal combined method achieves 9.50 days MAE on Indian monsoon onset (Nepal target)---a 32\% improvement over the rolling mean baseline and 9\% over CUSUM. Validated on the Lorenz system, logistic map, and three monsoon systems spanning both hemispheres (Indian/Nepal, Indian/Kerala, Western North Pacific), plus ENSO and a synthetic EEG dataset, the framework adds value precisely when the transition is gradual or obscured by noise.

\end{abstract}

\section{Introduction}\label{sec:introduction}

Bifurcations are the organizing centers of nonlinear dynamics. When a control parameter crosses a critical value, the qualitative structure of the attractor changes: a stable equilibrium gives way to oscillation, a periodic orbit destabilizes into chaos, or one chaotic attractor is replaced by another with different topology \citep{guckenheimer2013nonlinear, strogatz2015nonlinear}. These transitions govern phenomena across the sciences---from turbulence onset \citep{lorenz1963deterministic} and climate tipping points \citep{scheffer2009early} to epileptic seizures \citep{iasemidis2003adaptive} and cardiac arrhythmias \citep{glass2001synchronization}. Yet in most experimental and observational settings, the control parameter is unknown, the governing equations are unavailable, and the system is observed only through a noisy scalar time series. The fundamental question is: \emph{can we detect, from the time series alone, when the underlying attractor has changed?}

This paper develops a topological answer to that question. The core idea is simple: embed the time series in phase space via Takens delay coordinates, place a sliding window at a candidate transition time, and measure the \emph{topological complexity of the overlap} between the pre- and post-window point clouds as a function of scale. When both windows sample the same attractor, the overlap is topologically coherent; when a bifurcation intervenes, the overlap topology is disrupted, and a detection statistic spikes. We make this precise using the \emph{Mixup Euler Characteristic Profile} (Mixup ECP)---the Euler characteristic of the geometric intersection of ball unions around adjacent trajectory segments, viewed as a function of filtration scale. The Mixup ECP is defined and studied in \citep{majhi2026mixup_theory}; the present paper borrows only the definition and the Intersection Theorem, developing the rest---the detection pipeline, complementary signals, permutation inference, and multi-delay extension---from scratch.

The Mixup ECP is developed systematically in the forthcoming paper \emph{The Mixup Euler Characteristic Profile: Euler Calculus, Stability, and Integral Geometry for Topological Interaction of Ball Unions} \citep{majhi2026mixup_theory}, where it is shown to arise naturally from the Euler calculus on constructible functions, to satisfy $1$-Lipschitz interleaving stability under Hausdorff perturbations, and to extend to $k \geq 2$ point clouds via the product formula for the Euler integral. Geometric bounds on $|\Delta\chi|$ in terms of covering numbers and reach, Niyogi--Smale--Weinberger sampling guarantees, and an extension to Riemannian manifolds are also established there.

\subsection{Background and motivation}\label{sec:regime_problem}

A dynamical system undergoes a \emph{regime transition} when its qualitative behavior changes---in the language of bifurcation theory, when a parameter crosses a codimension-one bifurcation set: period-doubling cascades, Hopf bifurcations, crises, and boundary metamorphoses \citep{guckenheimer2013nonlinear, ott2002chaos}. In many settings, one observes only a scalar time series $\{x_t\}_{t=1}^N$ from a potentially high-dimensional system, and must infer \emph{when} (and whether) the underlying attractor has changed (Figure~\ref{fig:regime_problem}). This is the \textbf{regime detection problem}.

The problem is difficult for three reasons. First, the observation is low-dimensional: a scalar time series is a one-dimensional projection of the full state space. Second, real systems are noisy: observational noise, dynamical noise, and nonstationarity obscure the transition signal. Third, and most fundamentally, many regime transitions are \emph{gradual}: the system drifts through a bifurcation over an extended interval, so there is no single instant at which the regime ``switches.'' The challenge is not merely detecting that something changed, but localizing a transition that unfolds continuously in parameter space while manifesting as a qualitative change in attractor structure.

We take a moment to attempt to make precise what we mean by a \textit{regime} and a \textit{regime transition} (phase transition), as the terms are used somewhat loosely in the literature. At the most fundamental level, a regime is an attractor\textemdash{}an invariant set toward which nearby trajectories converge\textemdash{}together with its topological type: a fixed point, limit cycle, torus, or strange attractor. Two trajectory segments belong to the same regime if they are attracted to topologically equivalent sets, in the sense that their delay-embedded point clouds are homeomorphic at the relevant scales. A weaker, statistical notion identifies a regime with an ergodic component of the invariant measure: two states are in the same regime if their long-run time averages agree. Weakest of all is the phenomenological notion: a regime is a sustained, recognizable pattern in the observed time series\textemdash{}periodicity, quasiperiodicity, chaos. These three definitions are not equivalent; a chaotic attractor can contain multiple ergodic components, and two attractors can be topologically distinct yet phenomenologically indistinguishable in a scalar projection. A \textit{regime transition} is then a qualitative change at the topological level\textemdash{}the crossing of a codimension-one bifurcation set in parameter space\textemdash{}which may or may not manifest as a visible change at the statistical or phenomenological level. This hierarchy motivates the topological approach: statistical and signal-level methods are blind to transitions that do not alter univariate statistics, while attractor topology changes at every bifurcation by definition.

The term \textit{regime transition} as used in this paper is the observational counterpart of a bifurcation in the classical sense of Guckenheimer and Holmes \citep[Section~1.7]{guckenheimer2013nonlinear}. In that treatment, a bifurcation is defined precisely via structural stability: a vector field is structurally stable if all nearby vector fields are topologically equivalent to it, and a bifurcation occurs at parameter values where structural stability fails\textemdash{}that is, at the crossing of a codimension-one bifurcation set in parameter space, at which the topological orbit structure of the flow changes. Two vector fields on either side of the bifurcation set belong to different qualitative equivalence classes\textemdash{}they are not topologically conjugate\textemdash{}, and this inequivalence is the mathematical content of what we mean by different regimes. Strogatz \citep{strogatz2015nonlinear} captures the same idea more accessibly: bifurcations are qualitative changes in the dynamics as parameters vary, cataloged by the type of attractor created, destroyed, or transformed, without requiring the full machinery of structural stability. The difficulty is that these definitions are inaccessible in practice. They require knowledge of the governing equations, the parameter space, and the bifurcation set\textemdash{}none of which are available when the system is observed only through a noisy scalar time series. The term ``regime transition'' enters the applied sciences\textemdash{}meteorology, climate dynamics, neuroscience\textemdash{}precisely to fill this gap: it names, at the level of observation, what bifurcation theory names at the level of equations. A regime transition is what a bifurcation looks like from the outside. The Mixup Euler Characteristic Profile is designed for this observational setting: it provides a detection statistic sensitive to changes in attractor topology\textemdash{}the kind of change that underlies bifurcations in the sense of Guckenheimer and Holmes\textemdash{}without access to the governing equations, the parameter space, or the bifurcation set, working entirely from the reconstructed phase-space geometry of a finite time series. The Euler characteristic, while coarser than full topological conjugacy, is computationally tractable and stable and captures the essential qualitative distinction between regimes in which the bifurcation alters the number of connected components or cycles of the delay-embedded attractor.


The two preceding paragraphs suggest a natural question: can one give a TDA definition of regime transition directly in terms of the Mixup machinery? The Mixup barcode of \citet{wagner2024mixup}\textemdash{}the persistent homology of the intersection $U(X;r) \cap U(Y;r)$ across all scales, computed via inclusion--exclusion on Rips complexes\textemdash{}is a different invariant from the Mixup ECP: it encodes birth--death pairs of homological features rather than the Euler characteristic at each scale. This suggests the following definition: two adjacent trajectory segments $X$ and $Y$ belong to the \textit{same regime} if the persistent homology of $U(X;r) \cap U(Y;r)$ is close, in bottleneck distance, to the persistent homology of a single attractor sampled twice. A regime transition would then be a time 
$t$ at which this distance exceeds a threshold. Making this precise requires establishing a same-attractor limit for the Mixup barcode\textemdash{}a result analogous to the ECP limit $\Delta \chi(r) \to \chi(A)$ of the same-attractor property (Section~\ref{sec:definition})\textemdash{}and is a natural direction for future work. The Mixup ECP pursued in this paper requires only three Alpha complex computations per time step, versus a full persistent homology pipeline, and, in our experiments, classifies regimes at least as well as, and in some cases better than, the Mixup barcode, suggesting that the Euler characteristic captures the topologically relevant information for the transitions studied here. 

Moreover, the Mixup ECP produces a real-valued detection statistic $S(t)=\max_{r} |\Delta_\chi (r, X_t,Y_t)|$ directly, whereas converting a Mixup barcode to a scalar signal requires an additional vectorization step\textemdash{}such as taking a norm of the persistence diagram or computing bottleneck distance to a reference\textemdash{}and this vectorization inevitably discards information encoded in the full birth--death structure of the barcode. Finally, the Intersection Theorem (Theorem~\ref{thm:intersection}) gives the Mixup ECP a direct geometric grounding: $\Delta \chi(r; X, Y) = \chi\big(U(X;r) \cap U(Y;r)\big)$ is not an abstract algebraic summary, but the Euler characteristic of the actual geometric intersection of the two ball unions in $\mathbb{R}^d$, computed directly from the point clouds without any intermediate topological construction.

\subsection{Three notions of ``regime" and the role of the Mixup ECP}
We illustrate the three-level hierarchy and the role of the Mixup ECP with the logistic map $x_{n+1} = \lambda x_n(1-x_n)$, which is central to the validation experiments of Section~\ref{sec:logistic}. Consider three regimes: period-2 ($\lambda = 3.2$), period-4 ($\lambda = 3.5$), and chaos ($\lambda = 3.9$).

\textbf{Topological level.} The delay-embedded attractor changes type at each bifurcation. The period-2 orbit visits two distinct locations repeatedly, so the delay-embedded point cloud concentrates into two clusters\textemdash{}topologically, two disjoint components. The period-4 orbit produces four such clusters. The chaotic attractor fills a connected region of phase space, giving a single connected component. No two of these are homeomorphic: two disjoint components cannot be continuously deformed into four, and no disconnected set is homeomorphic to a connected one. The topological definition therefore correctly identifies all three as distinct regimes, and identifies a regime transition at each bifurcation.

\textbf{Statistical level.} The invariant measures of the period-2 and period-4 orbits are supported on different sets, so their ergodic components technically differ. However, both orbits are confined to roughly the same region of phase space, with nearly identical means and variances ($\bar{x} \approx 1.4$, $\sigma^2 \approx 0.5$ for both; see Figure~\ref{fig:why_topology}(a)). A statistical test based on moments would fail to reject the null that the two segments come from the same distribution. The period-2 to period-4 transition is thus invisible at the statistical level in practice, even though it is present in principle. The period-4 to chaos transition is more visible statistically, as the chaotic orbit fills a broader region of phase space, but detection is unreliable near the bifurcation point where the chaotic band is still narrow.

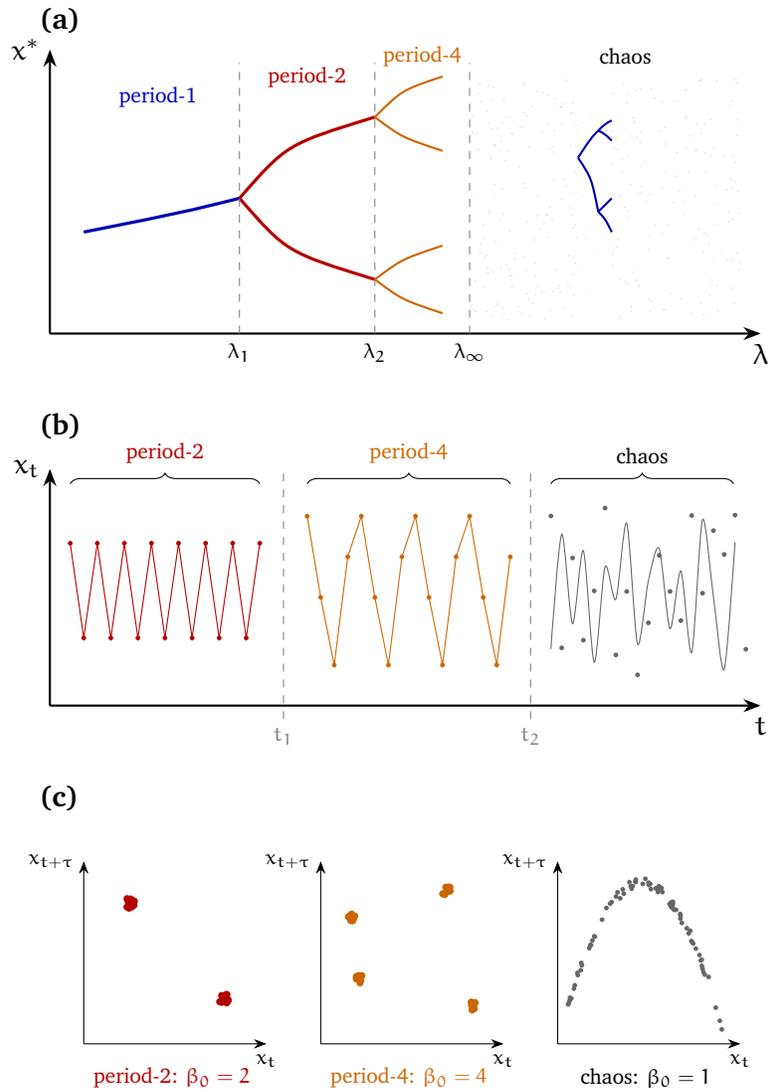
\begin{figure}[htbp]
\centering
\begin{tikzpicture}[>=Stealth,scale=0.9]

\begin{scope}[shift={(0,0)}]
  \node[anchor=south west,font=\bfseries] at (-0.3,4.3) {(a)};
  \draw[->,thick] (0,0) -- (10.5,0) node[below] {$\lambda$};
  \draw[->,thick] (0,0) -- (0,4.2) node[left] {$x^*$};

  \draw[blue!70!black, very thick] (0.5,1.5) .. controls (2,1.8) .. (2.8,2.0);

  \draw[red!70!black, very thick] (2.8,2.0) .. controls (3.5,2.8) .. (4.8,3.2);
  \draw[red!70!black, very thick] (2.8,2.0) .. controls (3.5,1.2) .. (4.8,0.8);

  \draw[orange!80!black, thick] (4.8,3.2) .. controls (5.2,3.6) .. (5.8,3.8);
  \draw[orange!80!black, thick] (4.8,3.2) .. controls (5.2,2.9) .. (5.8,2.7);
  \draw[orange!80!black, thick] (4.8,0.8) .. controls (5.2,1.1) .. (5.8,1.3);
  \draw[orange!80!black, thick] (4.8,0.8) .. controls (5.2,0.5) .. (5.8,0.3);

  \pgfmathsetseed{42}
  \foreach \i in {1,...,400} {
    \pgfmathsetmacro{\lam}{6.2 + 4.0*rnd}
    \pgfmathsetmacro{\yy}{0.2 + 3.6*rnd}
    \pgfmathsetmacro{\prob}{rnd}
    \ifnum\i<350
      \fill[black!50, opacity=0.3] (\lam,\yy) circle (0.3pt);
    \fi
  }

  \fill[white,opacity=0.7] (7.8,0.1) rectangle (8.3,4.0);
  \draw[blue!70!black, thick] (7.8,2.6) .. controls (8.0,2.9) .. (8.1,3.0);
  \draw[blue!70!black, thick] (7.8,2.6) .. controls (8.0,2.3) .. (8.1,1.8);
  \draw[blue!70!black, thick] (8.1,3.0) .. controls (8.2,3.1) .. (8.3,3.15);
  \draw[blue!70!black, thick] (8.1,3.0) .. controls (8.2,2.95) .. (8.3,2.85);
  \draw[blue!70!black, thick] (8.1,1.8) .. controls (8.2,1.9) .. (8.3,2.0);
  \draw[blue!70!black, thick] (8.1,1.8) .. controls (8.2,1.7) .. (8.3,1.5);

  \draw[dashed, gray] (2.8,0) -- (2.8,4.0);
  \draw[dashed, gray] (4.8,0) -- (4.8,4.0);
  \draw[dashed, gray] (6.2,0) -- (6.2,4.0);
  \node[below,font=\scriptsize] at (2.8,0) {$\lambda_1$};
  \node[below,font=\scriptsize] at (4.8,0) {$\lambda_2$};
  \node[below,font=\scriptsize] at (6.2,0) {$\lambda_\infty$};

  \node[font=\scriptsize,blue!70!black] at (1.6,3.5) {period-1};
  \node[font=\scriptsize,red!70!black] at (3.8,3.8) {period-2};
  \node[font=\scriptsize,orange!80!black] at (5.5,4.1) {period-4};
  \node[font=\scriptsize] at (8.5,4.1) {chaos};
\end{scope}

\begin{scope}[shift={(0,-5.5)}]
  \node[anchor=south west,font=\bfseries] at (-0.3,3.8) {(b)};
  \draw[->,thick] (0,0) -- (10.5,0) node[below] {$t$};
  \draw[->,thick] (0,0) -- (0,3.5) node[left] {$x_t$};

  \foreach \i in {0,...,14} {
    \pgfmathsetmacro{\x}{0.3 + \i*0.2}
    \pgfmathsetmacro{\ymod}{mod(\i,2)}
    \pgfmathsetmacro{\y}{ifthenelse(\ymod==0, 2.4, 1.0)}
    \fill[red!70!black] (\x,\y) circle (1pt);
  }
  \draw[red!70!black,thin] (0.3,2.4) -- (0.5,1.0) -- (0.7,2.4) -- (0.9,1.0) -- (1.1,2.4) -- (1.3,1.0) -- (1.5,2.4) -- (1.7,1.0) -- (1.9,2.4) -- (2.1,1.0) -- (2.3,2.4) -- (2.5,1.0) -- (2.7,2.4) -- (2.9,1.0) -- (3.1,2.4);

  \foreach \i in {0,...,15} {
    \pgfmathsetmacro{\x}{3.8 + \i*0.2}
    \pgfmathsetmacro{\imod}{mod(\i,4)}
    \pgfmathsetmacro{\y}{ifthenelse(\imod==0, 2.8, ifthenelse(\imod==1, 1.6, ifthenelse(\imod==2, 0.6, 2.2)))}
    \fill[orange!80!black] (\x,\y) circle (1pt);
  }
  \draw[orange!80!black,thin] (3.8,2.8) -- (4.0,1.6) -- (4.2,0.6) -- (4.4,2.2) -- (4.6,2.8) -- (4.8,1.6) -- (5.0,0.6) -- (5.2,2.2) -- (5.4,2.8) -- (5.6,1.6) -- (5.8,0.6) -- (6.0,2.2) -- (6.2,2.8) -- (6.4,1.6) -- (6.6,0.6) -- (6.8,2.2);

  \pgfmathsetseed{99}
  \def\prevx{7.4}\def\prevy{1.8}
  \foreach \i in {0,...,18} {
    \pgfmathsetmacro{\x}{7.4 + \i*0.16}
    \pgfmathsetmacro{\y}{0.4 + 2.6*rnd}
    \fill[black!60] (\x,\y) circle (1pt);
    \xdef\currx{\x}\xdef\curry{\y}
  }
  \pgfmathsetseed{99}
  \draw[black!60,thin] plot[smooth] coordinates {
    (7.40,{0.4+2.6*0.168}) (7.56,{0.4+2.6*0.821}) (7.72,{0.4+2.6*0.312})
    (7.88,{0.4+2.6*0.727}) (8.04,{0.4+2.6*0.092}) (8.20,{0.4+2.6*0.634})
    (8.36,{0.4+2.6*0.451}) (8.52,{0.4+2.6*0.883}) (8.68,{0.4+2.6*0.197})
    (8.84,{0.4+2.6*0.556}) (9.00,{0.4+2.6*0.738}) (9.16,{0.4+2.6*0.284})
    (9.32,{0.4+2.6*0.611}) (9.48,{0.4+2.6*0.149}) (9.64,{0.4+2.6*0.903})
    (9.80,{0.4+2.6*0.467}) (9.96,{0.4+2.6*0.052}) (10.12,{0.4+2.6*0.775})
  };

  \draw[decorate,decoration={brace,amplitude=4pt,raise=2pt}] (0.3,3.2) -- (3.1,3.2)
    node[midway,above=6pt,font=\scriptsize,red!70!black] {period-2};
  \draw[decorate,decoration={brace,amplitude=4pt,raise=2pt}] (3.8,3.2) -- (6.8,3.2)
    node[midway,above=6pt,font=\scriptsize,orange!80!black] {period-4};
  \draw[decorate,decoration={brace,amplitude=4pt,raise=2pt}] (7.4,3.2) -- (10.1,3.2)
    node[midway,above=6pt,font=\scriptsize] {chaos};

  \draw[dashed,gray] (3.45,-0.2) -- (3.45,3.5);
  \draw[dashed,gray] (7.1,-0.2) -- (7.1,3.5);
  \node[font=\scriptsize,gray] at (3.45,-0.45) {$t_1$};
  \node[font=\scriptsize,gray] at (7.1,-0.45) {$t_2$};
\end{scope}

\begin{scope}[shift={(0,-10.5)}]
  \node[anchor=south west,font=\bfseries] at (-0.3,3.3) {(c)};

  \begin{scope}[shift={(0.5,0)}]
    \draw[->,thin] (0,0) -- (2.7,0) node[below,font=\scriptsize] {$x_t$};
    \draw[->,thin] (0,0) -- (0,2.7) node[left,font=\scriptsize] {$x_{t+\tau}$};
    \foreach \i in {1,...,12} {
      \pgfmathsetmacro{\ax}{0.6+0.15*rnd}
      \pgfmathsetmacro{\ay}{2.0+0.15*rnd}
      \fill[red!70!black] (\ax,\ay) circle (1.5pt);
    }
    \foreach \i in {1,...,12} {
      \pgfmathsetmacro{\bx}{2.0+0.15*rnd}
      \pgfmathsetmacro{\by}{0.6+0.15*rnd}
      \fill[red!70!black] (\bx,\by) circle (1.5pt);
    }
    \node[font=\scriptsize,red!70!black] at (1.3,-0.5) {period-2: $\beta_0 = 2$};
  \end{scope}

  \begin{scope}[shift={(4.0,0)}]
    \draw[->,thin] (0,0) -- (2.7,0) node[below,font=\scriptsize] {$x_t$};
    \draw[->,thin] (0,0) -- (0,2.7) node[left,font=\scriptsize] {$x_{t+\tau}$};
    \foreach \i in {1,...,8} {
      \pgfmathsetmacro{\ax}{0.4+0.12*rnd}
      \pgfmathsetmacro{\ay}{1.8+0.12*rnd}
      \fill[orange!80!black] (\ax,\ay) circle (1.5pt);
    }
    \foreach \i in {1,...,8} {
      \pgfmathsetmacro{\bx}{1.8+0.12*rnd}
      \pgfmathsetmacro{\by}{2.2+0.12*rnd}
      \fill[orange!80!black] (\bx,\by) circle (1.5pt);
    }
    \foreach \i in {1,...,8} {
      \pgfmathsetmacro{\cx}{2.2+0.12*rnd}
      \pgfmathsetmacro{\cy}{0.5+0.12*rnd}
      \fill[orange!80!black] (\cx,\cy) circle (1.5pt);
    }
    \foreach \i in {1,...,8} {
      \pgfmathsetmacro{\dx}{0.5+0.12*rnd}
      \pgfmathsetmacro{\dy}{0.9+0.12*rnd}
      \fill[orange!80!black] (\dx,\dy) circle (1.5pt);
    }
    \node[font=\scriptsize,orange!80!black] at (1.3,-0.5) {period-4: $\beta_0 = 4$};
  \end{scope}

  \begin{scope}[shift={(7.5,0)}]
    \draw[->,thin] (0,0) -- (2.7,0) node[below,font=\scriptsize] {$x_t$};
    \draw[->,thin] (0,0) -- (0,2.7) node[left,font=\scriptsize] {$x_{t+\tau}$};
    \pgfmathsetseed{77}
    \foreach \i in {1,...,80} {
      \pgfmathsetmacro{\xx}{0.15 + 2.3*rnd}
      \pgfmathsetmacro{\xnorm}{\xx/2.5}
      \pgfmathsetmacro{\ybase}{3.8*\xnorm*(1-\xnorm)*2.5}
      \pgfmathsetmacro{\yy}{\ybase + 0.15*(rnd-0.5)}
      \fill[black!60] (\xx,\yy) circle (1pt);
    }
    \node[font=\scriptsize] at (1.3,-0.5) {chaos: $\beta_0 = 1$};
  \end{scope}
\end{scope}

\end{tikzpicture}
\caption{The regime detection problem, illustrated on the logistic map $x_{n+1} = \lambda x_n(1-x_n)$. (a)~Bifurcation diagram: the attractor structure changes qualitatively as $\lambda$ increases through period-doubling cascades into chaos. (b)~Observed time series when $\lambda$ drifts: the scalar signal shows increasing complexity, but pinpointing transitions from the time series alone is difficult. (c)~Delay-embedded point clouds ($\tau=1$, $d=2$) from three regimes: the attractor \emph{topology} changes---from two clusters (period-2) to four clusters (period-4) to a connected chaotic band---even when signal amplitude remains comparable.}
\label{fig:regime_problem}
\end{figure}

\textbf{Phenomenological level.} The scalar time series in all three regimes consists of bounded oscillations of comparable amplitude. The doubling of the period from 2 to 4 is subtle and easily masked by noise; a practitioner inspecting the raw signal would likely see no transition. The transition to chaos produces increasing irregularity, but pinpointing the bifurcation from the scalar signal alone is difficult, as Figure~\ref{fig:regime_problem}(b) illustrates.

\textbf{Mixup ECP.} The Euler characteristic of the delay-embedded point cloud takes the value $\chi = \beta_0 - \beta_1 = 2$, $4$, and $1$ for period-2, period-4, and chaos, respectively. When the Mixup ECP is computed between adjacent sliding windows straddling a bifurcation\textemdash{}one window in the period-2 regime, one in the period-4 regime\textemdash{}the ball unions around the two clusters of the period-2 cloud and the four clusters of the period-4 cloud overlap weakly, producing a small $\Delta\chi(r)$ profile and a low detection statistic $S(t) = \max_r |\Delta\chi(r)|$. By contrast, when both windows sample the same attractor, the ball unions overlap richly, and $S(t)$ is large. The permutation test of Section~\ref{sec:hypothesis} confirms this: same-regime comparisons (period-2 vs.\ period-2) produce $S_{\mathrm{obs}} = 6$ with null $6.0 \pm 0.0$ and $p = 1.0$, while cross-regime comparisons (period-2 vs.\ chaos) produce $S_{\mathrm{obs}} = 3$ with null $17.7 \pm 2.2$ and $p < 0.01$. The Euler characteristic thus captures the $\beta_0$ change at the period-doubling bifurcation cleanly. It does not, however, distinguish two chaotic attractors with identical Betti numbers but different invariant measures\textemdash{}the precise sense in which it is coarser than full topological conjugacy. In such cases, it is the full shape of the profile $r \mapsto \Delta\chi(r)$, rather than the scalar maximum $S$, that carries the regime information.

\subsection{Limitations of existing approaches to regime detection}\label{sec:existing}

Several families of methods address regime detection, each capturing a different aspect of the dynamics.

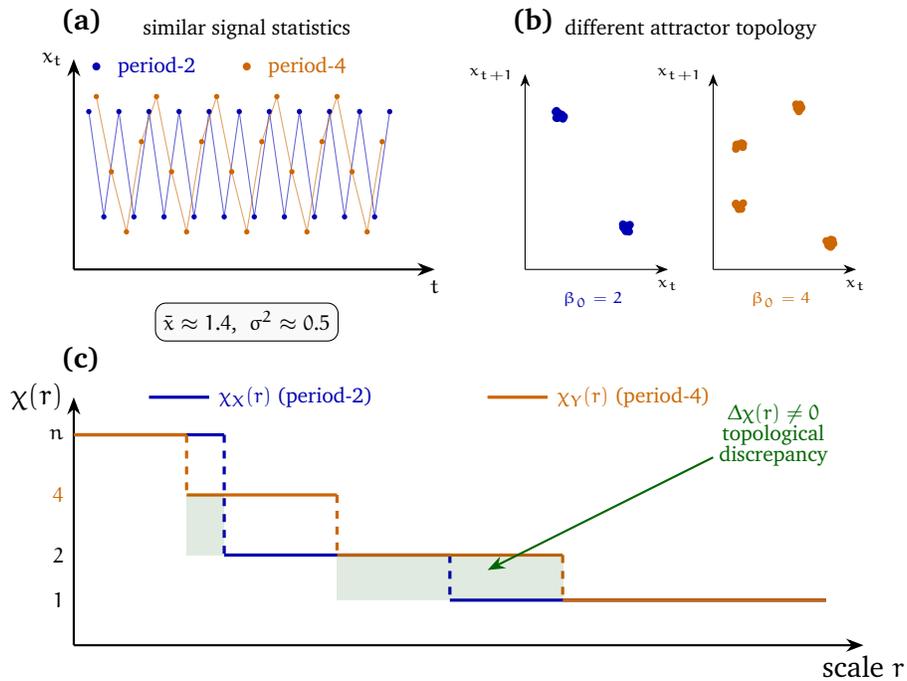
\begin{figure}[!htbp]
\centering
\begin{tikzpicture}[>=Stealth]
\begin{scope}[shift={(0,0)}]
  \node[anchor=south west,font=\bfseries] at (-0.3,3.0) {(a)};
  \draw[->,thick] (0,0) -- (4.8,0) node[below,font=\scriptsize] {$t$};
  \draw[->,thick] (0,0) -- (0,2.8) node[left,font=\scriptsize] {$x_t$};
  \node[font=\scriptsize] at (2.3,3.2) {similar signal statistics};

  \foreach \i in {0,...,20} {
    \pgfmathsetmacro{\x}{0.2 + \i*0.2}
    \pgfmathsetmacro{\imod}{mod(\i,2)}
    \pgfmathsetmacro{\y}{ifthenelse(\imod==0, 2.1, 0.7)}
    \fill[blue!70!black] (\x,\y) circle (1pt);
  }
  \draw[blue!70!black,thin,opacity=0.6] (0.2,2.1)--(0.4,0.7)--(0.6,2.1)--(0.8,0.7)--(1.0,2.1)--(1.2,0.7)--(1.4,2.1)--(1.6,0.7)--(1.8,2.1)--(2.0,0.7)--(2.2,2.1)--(2.4,0.7)--(2.6,2.1)--(2.8,0.7)--(3.0,2.1)--(3.2,0.7)--(3.4,2.1)--(3.6,0.7)--(3.8,2.1)--(4.0,0.7)--(4.2,2.1);

  \foreach \i in {0,...,19} {
    \pgfmathsetmacro{\x}{0.3 + \i*0.2}
    \pgfmathsetmacro{\imod}{mod(\i,4)}
    \pgfmathsetmacro{\y}{ifthenelse(\imod==0, 2.3, ifthenelse(\imod==1, 1.3, ifthenelse(\imod==2, 0.5, 1.7)))}
    \fill[orange!80!black] (\x,\y) circle (1pt);
  }
  \draw[orange!80!black,thin,opacity=0.6] (0.3,2.3)--(0.5,1.3)--(0.7,0.5)--(0.9,1.7)--(1.1,2.3)--(1.3,1.3)--(1.5,0.5)--(1.7,1.7)--(1.9,2.3)--(2.1,1.3)--(2.3,0.5)--(2.5,1.7)--(2.7,2.3)--(2.9,1.3)--(3.1,0.5)--(3.3,1.7)--(3.5,2.3)--(3.7,1.3)--(3.9,0.5)--(4.1,1.7);

  \node[draw,rounded corners,font=\scriptsize,inner sep=3pt,fill=gray!5] at (2.3,-0.7) {%
    $\bar{x} \approx 1.4$,\; $\sigma^2 \approx 0.5$%
  };

  \fill[blue!70!black] (0.3,2.7) circle (1.5pt);
  \node[right,font=\scriptsize,blue!70!black] at (0.45,2.7) {period-2};
  \fill[orange!80!black] (2.3,2.7) circle (1.5pt);
  \node[right,font=\scriptsize,orange!80!black] at (2.45,2.7) {period-4};
\end{scope}

\begin{scope}[shift={(6,0)}]
  \node[anchor=south west,font=\bfseries] at (-0.3,3.0) {(b)};
  \node[font=\scriptsize] at (2.2,3.2) {different attractor topology};

  \begin{scope}[shift={(0,0)}]
    \draw[->,thin] (0,0) -- (1.9,0) node[below,font=\tiny] {$x_t$};
    \draw[->,thin] (0,0) -- (0,2.6) node[left,font=\tiny] {$x_{t+1}$};
    \pgfmathsetseed{11}
    \foreach \i in {1,...,15} {
      \fill[blue!70!black] ({0.4+0.12*rnd},{2.0+0.12*rnd}) circle (1.5pt);
      \fill[blue!70!black] ({1.3+0.12*rnd},{0.5+0.12*rnd}) circle (1.5pt);
    }
    \node[font=\tiny,blue!70!black] at (0.9,-0.4) {$\beta_0=2$};
  \end{scope}

  \begin{scope}[shift={(2.5,0)}]
    \draw[->,thin] (0,0) -- (1.9,0) node[below,font=\tiny] {$x_t$};
    \draw[->,thin] (0,0) -- (0,2.6) node[left,font=\tiny] {$x_{t+1}$};
    \pgfmathsetseed{22}
    \foreach \i in {1,...,10} {
      \fill[orange!80!black] ({0.3+0.1*rnd},{1.6+0.1*rnd}) circle (1.5pt);
      \fill[orange!80!black] ({1.1+0.1*rnd},{2.1+0.1*rnd}) circle (1.5pt);
      \fill[orange!80!black] ({1.5+0.1*rnd},{0.3+0.1*rnd}) circle (1.5pt);
      \fill[orange!80!black] ({0.3+0.1*rnd},{0.8+0.1*rnd}) circle (1.5pt);
    }
    \node[font=\tiny,orange!80!black] at (0.9,-0.4) {$\beta_0=4$};
  \end{scope}
\end{scope}

\begin{scope}[shift={(0,-5.0)}]
  \node[anchor=south west,font=\bfseries] at (-0.3,3.5) {(c)};
  \draw[->,thick] (0,0) -- (10.5,0) node[below] {scale $r$};
  \draw[->,thick] (0,0) -- (0,3.3) node[left] {$\chi(r)$};

  \draw[blue!70!black, very thick]
    (0,2.8) -- (2.0,2.8)    
    (2.0,1.2) -- (5.0,1.2)  
    (5.0,0.6) -- (10,0.6);  
  \draw[blue!70!black, very thick, dashed] (2.0,2.8) -- (2.0,1.2);
  \draw[blue!70!black, very thick, dashed] (5.0,1.2) -- (5.0,0.6);

  \draw[orange!80!black, very thick]
    (0,2.8) -- (1.5,2.8)
    (1.5,2.0) -- (3.5,2.0)  
    (3.5,1.2) -- (6.5,1.2)  
    (6.5,0.6) -- (10,0.6);  
  \draw[orange!80!black, very thick, dashed] (1.5,2.8) -- (1.5,2.0);
  \draw[orange!80!black, very thick, dashed] (3.5,2.0) -- (3.5,1.2);
  \draw[orange!80!black, very thick, dashed] (6.5,1.2) -- (6.5,0.6);

  \node[left,font=\scriptsize] at (0,2.8) {$n$};
  \node[left,font=\scriptsize,orange!80!black] at (0,2.0) {$4$};
  \node[left,font=\scriptsize] at (0,1.2) {$2$};
  \node[left,font=\scriptsize] at (0,0.6) {$1$};

  \fill[green!30!black, opacity=0.12]
    (1.5,2.0) -- (2.0,2.0) -- (2.0,1.2) -- (1.5,1.2) -- cycle;
  \fill[green!30!black, opacity=0.12]
    (3.5,1.2) -- (5.0,1.2) -- (5.0,0.6) -- (3.5,0.6) -- cycle;
  \fill[green!30!black, opacity=0.12]
    (5.0,1.2) -- (6.5,1.2) -- (6.5,0.6) -- (5.0,0.6) -- cycle;

  \draw[->,thick,green!40!black] (8.5,2.5) -- (5.5,1.0);
  \node[font=\scriptsize,green!40!black,align=center] at (9.3,2.8) {$\Delta\chi(r) \neq 0$\\[-2pt]\scriptsize topological\\[-2pt]\scriptsize discrepancy};

  \draw[blue!70!black, very thick] (1.0,3.3) -- (1.8,3.3);
  \node[right,font=\scriptsize,blue!70!black] at (1.8,3.3) {$\chi_X(r)$ (period-2)};
  \draw[orange!80!black, very thick] (5.5,3.3) -- (6.3,3.3);
  \node[right,font=\scriptsize,orange!80!black] at (6.3,3.3) {$\chi_Y(r)$ (period-4)};
\end{scope}

\end{tikzpicture}
\caption{Why signal-level methods miss topological transitions. (a)~Two segments of the logistic map---period-2 ($\lambda = 3.3$, blue) and period-4 ($\lambda = 3.5$, orange)---have nearly identical mean and variance; signal-level change-point methods see no transition. (b)~Delay embeddings reveal different attractor topologies: two clusters ($\beta_0 = 2$) vs.\ four clusters ($\beta_0 = 4$). (c)~The Euler characteristic surfaces differ: both start at $\chi = n$ and end at $\chi = 1$, but the period-4 ECS passes through an intermediate $\chi = 4$ plateau. The shaded regions mark scales where $\Delta\chi(r) \neq 0$---the Mixup ECS captures this topological discrepancy.}
\label{fig:why_topology}
\end{figure}

\textbf{Signal-level methods.} Classical change-point detection---CUSUM \citep{page1954continuous}, rolling mean shifts, maximum slope methods---operates on the raw time series or simple derived statistics (mean, variance, spectral power). These methods are effective when the transition produces a clear amplitude or frequency shift, but are blind to changes in attractor geometry that do not alter univariate signal statistics (Figure~\ref{fig:why_topology}).

\textbf{Early warning signals.} The critical slowing down framework \citep{scheffer2009early, dakos2008slowing} detects increased autocorrelation and variance as a system approaches a tipping point. This approach is grounded in bifurcation theory (specifically, the slowing of the dominant eigenvalue near a fold bifurcation) but assumes a specific bifurcation mechanism and requires long pre-transition baselines. It does not apply to transitions without critical slowing down---for instance, period-doubling cascades or interior crises.

\textbf{Recurrence-based methods.} Recurrence quantification analysis (RQA) \citep{marwan2007recurrence, eckmann1995recurrence} and recurrence networks characterize the structure of phase-space recurrences. These methods detect regime transitions through changes in recurrence statistics (determinism, laminarity, trapping time) and have been applied to climate \citep{marwan2002nonlinear} and physiological signals. However, they rely on a single distance threshold, collapsing the multiscale structure of the attractor into scalar statistics.

\textbf{Lyapunov exponents and dimension estimates.} Finite-time Lyapunov exponents \citep{wolf1985determining, rosenstein1993practical} and correlation dimension estimates \citep{grassberger1983characterization} track the divergence rate and fractal structure of the attractor. These are gold-standard dynamical diagnostics, but they require long, stationary trajectory segments for reliable estimation---precisely the condition violated at a regime transition.

A common limitation of all these approaches is that they track \emph{statistical} or \emph{metrical} properties of the dynamics (means, variances, divergence rates, recurrence frequencies) rather than the \emph{shape} of the attractor. Yet bifurcations are fundamentally topological events: a period-doubling bifurcation creates a new connected component in Poincar\'e section; a Hopf bifurcation creates a limit cycle from an equilibrium; a crisis merges or destroys an attractor. The attractor's topology changes, but existing detection methods access this only indirectly.

\textbf{Topological data analysis.} Persistent homology applied to delay-embedded time series offers a principled way to track attractor shape across scales. The theoretical foundations for this program were laid by \citet{perea2015sliding}, who established the connection between sliding-window embeddings and persistent homology, with convergence guarantees relating maximum persistence to signal periodicity (see also \citet{perea2019topological} for a survey). In the bifurcation detection setting, Khasawneh, Munch, and collaborators have developed a systematic program: BuZZ uses zigzag persistent homology to detect Hopf bifurcations in a single persistence computation \citep{tymochko2020using}; CROCKER plots visualize how Betti numbers evolve across a one-parameter family of Rips filtrations, enabling data-driven bifurcation detection across ten dynamical systems \citep{guzel2022detecting}; and the framework has been extended to stochastic systems, where homological bifurcation plots detect phenomenological (P-type) bifurcations via superlevel persistence of estimated densities \citep{tanweer2024topological, tanweer2024unreliable}. In climate science, \citet{strommen2023topological} equate weather regimes with nontrivial topological structure of the attractor and use persistent homology to detect regimes in the Lorenz system and ERA5 reanalysis data; \citet{faranda2024climate} survey the broader interface between TDA and dynamical systems for climate applications.

Parallel to the persistent homology program, the Euler characteristic has emerged as a computationally efficient topological summary. The Euler Characteristic Transform \citep{turner2014persistent} is injective on shapes; Euler characteristic curves and profiles provide stable, distributable alternatives to persistence for large-scale problems \citep{dlotko2023euler, hacquard2024euler}; and Euler characteristic surfaces have been applied to characterize fluid dynamical systems \citep{roy2023characterizing} and classify time series \citep{luwang2026interpretable}.

These methods represent significant advances, but they share a structural limitation: they characterize the topology of a \emph{single} point cloud (one window of the time series) and detect transitions by comparing successive single-cloud summaries. Our approach differs architecturally. The Mixup ECP computes the topology of the \emph{interaction} between two adjacent windows---specifically, the Euler characteristic of the geometric intersection $\mathcal{U}(X;r) \cap \mathcal{U}(Y;r)$ (Theorem~\ref{thm:intersection}). This two-sample construction has a built-in null: the detection statistic is zero when the ball unions do not overlap (dead-zone property), and converges to $\chi(\mathcal{A})$ when both windows sample the same attractor. By contrast, single-cloud methods must compare topological summaries across time via an external distance or visualization (e.g., bottleneck distance between persistence diagrams, or visual inspection of CROCKER contours). The Mixup construction internalizes this comparison, producing a single detection statistic $S(t)$ with a natural permutation test. Computationally, the Mixup ECP on the diagonal requires three Alpha complex computations per time step---$O(n^2)$ in low-dimensional embeddings---compared to the $O(n^3)$ cost of full persistent homology pipelines.

\subsection{Topological approach to regime detection}\label{sec:topo_approach}

Topological data analysis (TDA) offers a complementary perspective: rather than tracking signal amplitude or divergence rates, TDA characterizes the \emph{shape} of data across scales \citep{edelsbrunner2010computational, perea2019topological}. Applied to time series via Takens delay embedding \citep{takens2006detecting, sauer1991embedology}, TDA can detect changes in the geometric and topological structure of reconstructed phase-space attractors---precisely the kind of changes that bifurcations produce. Recent work has demonstrated the value of topological methods for climate dynamics \citep{faranda2024climate, strommen2023topological}, flow regime classification \citep{koenig2026topological}, and time series analysis \citep{perea2019topological, luwang2026interpretable}. The present work builds on \citet{alvarado2025detecting}, who first demonstrated that persistent homology applied to delay-embedded monsoon index data can detect the Indian monsoon onset as a topological regime transition; the Mixup ECP framework developed here provides a computationally cheaper and theoretically grounded alternative to the persistent homology approach of that earlier study.

The \emph{Mixup Euler Characteristic Profile} (Mixup ECP) measures the topological interaction between point clouds as a function of scale. For two disjoint point clouds $X, Y \subset \R^d$, the Mixup ECP at scale $r$ is
\begin{equation}\label{eq:delta_chi_intro}
    \Delta\chi(r;\, X, Y) = \chi_X(r) + \chi_Y(r) - \chi_{X \cup Y}(r),
\end{equation}
where $\chi_X(r) = \chi(\mathcal{U}(X; r))$ is the Euler characteristic of the ball union $\mathcal{U}(X; r) = \bigcup_{p \in X} B(p, r)$ \citep{roy2025euler}, and similarly for $Y$ and $X \cup Y$. The Intersection Theorem (Theorem~\ref{thm:intersection}; see \citep{majhi2026mixup_theory} for proof) establishes that $\Delta\chi(r) = \chi(\mathcal{U}(X; r) \cap \mathcal{U}(Y; r))$---the Euler characteristic of the geometric intersection of ball unions. This identity gives the detection statistic $S(t) = \max_{r} |\Delta\chi(r)|$ a direct geometric meaning: it measures the topological complexity of the overlap region between adjacent trajectory segments. When the attractor changes across the window boundary---as it does at a bifurcation---the overlap topology is disrupted and $S(t)$ spikes.

\begin{remark}[Independent scales]
More generally, the Mixup ECP can be evaluated at independent scales $(r, s)$ for the two clouds, giving $\Delta\chi(r, s) = \chi(\mathcal{U}(X; r) \cap \mathcal{U}(Y; s))$ \citep{majhi2026mixup_theory}. This biparameter extension is natural when the pre- and post-transition attractors have different characteristic sizes, but requires weighted Alpha complexes. Throughout this paper we use the diagonal $r = s$, which requires only three standard Alpha complex computations per time step.
\end{remark}

The Mixup ECP is complemented by a geometric detection signal---the \emph{Complexity Variance} $G(t) = dV/dt$, where $V(t)$ is the variance of the sliding-window point cloud. We establish a formal connection: the variance controls the \emph{support width} of the ECP (Proposition~\ref{prop:var_ecs}), while the Mixup ECP captures the \emph{topological content} of the cross-window interaction. These signals are complementary by construction: a size-preserving topology change (e.g., period-doubling at fixed amplitude) is detected by $S$ but not $G$; a topology-preserving size change (e.g., amplitude growth on a fixed attractor) is detected by $G$ but not $S$.

\subsection{Contributions and outline}\label{sec:contributions}

This paper develops the Mixup ECP into a practical regime detection framework, and tests it systematically on systems of increasing complexity.

\begin{enumerate}
    \item \textbf{Detection pipeline.} We combine the Mixup ECP with three complementary signals---Complexity Variance, Higuchi fractal dimension, and rolling mean---into a four-signal detection framework with a formal connection between the topological and geometric signals through the ECP support (Section~\ref{sec:pipeline}).
    \item \textbf{Permutation test for regime detection.} We formalize the detection problem as a permutation test with low-side rejection: a regime transition is detected when the cross-window overlap topology is unusually weak. We establish validity, consistency, and validate on synthetic bifurcations (Lorenz, logistic map) where the test achieves perfect discrimination (Section~\ref{sec:hypothesis}).
    \item \textbf{Multi-delay extension.} Embedding at $m$ delays produces $2m$ point clouds whose Mixup ECP automatically selects the dynamical timescale at which the bifurcation is most topologically visible, provably dominating any fixed-delay test in power (Section~\ref{sec:multi_delay}).
    \item \textbf{Validation on canonical bifurcations.} We test on the Lorenz system (fixed point to chaos) and the logistic map (period-doubling cascade), including Lyapunov exponent comparison ($r = 0.68$) and noise robustness experiments (Section~\ref{sec:synthetic}).
    \item \textbf{Real-world regime transitions.} Monsoon onset detection across three systems (Indian/Nepal, Indian/Kerala, Western North Pacific)---topology helps on gradual-onset targets (32\% improvement on Nepal over rolling mean, 50\% on WNP) and is harmless on sharp-onset targets---plus ENSO phase transitions and synthetic EEG seizure onset (Section~\ref{sec:applications})..
\end{enumerate}

\textbf{Outline.} Section~\ref{sec:pipeline} describes the detection pipeline. Section~\ref{sec:framework} develops the theoretical foundations: properties of the Mixup ECP, the variance--ECP connection, the permutation test, and the multi-delay extension. Section~\ref{sec:synthetic} validates on synthetic dynamical systems. Section~\ref{sec:applications} demonstrates on geophysical and neural applications. Section~\ref{sec:discussion} discusses when topology helps. Section~\ref{sec:conclusion} concludes.

\section{Detection Pipeline}\label{sec:pipeline}

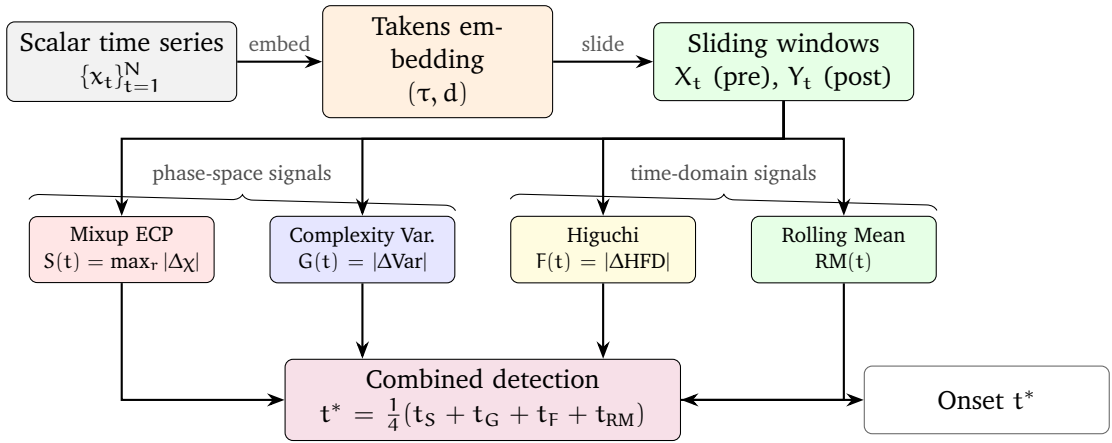
\begin{figure}[!htbp]
\centering
\begin{tikzpicture}[
    >=Stealth,
    box/.style={draw, rounded corners=3pt, minimum height=0.9cm, text width=2.8cm, align=center, font=\small},
    signal/.style={draw, rounded corners=3pt, minimum height=0.7cm, text width=2.2cm, align=center, font=\scriptsize, fill=blue!8},
    arrow/.style={->, thick},
    label/.style={font=\scriptsize, text=gray!70!black},
]

\node[box, fill=gray!10] (ts) at (0,0) {Scalar time series\\$\{x_t\}_{t=1}^N$};

\node[box, fill=orange!12] (emb) at (4.2,0) {Takens embedding\\$(\tau, d)$};

\node[box, fill=green!10, text width=3.2cm] (win) at (8.8,0) {Sliding windows\\$X_t$ (pre), $Y_t$ (post)};

\node[signal, fill=red!10] (S) at (0.0, -2.5) {Mixup ECP\\$S(t) = \max_r |\Delta\chi|$};
\node[signal, fill=blue!10] (G) at (3.2, -2.5) {Complexity Var.\\$G(t) = |\Delta\mathrm{Var}|$};
\node[signal, fill=yellow!15] (F) at (6.4, -2.5) {Higuchi\\$F(t) = |\Delta\mathrm{HFD}|$};
\node[signal, fill=green!10] (RM) at (9.6, -2.5) {Rolling Mean\\$\mathrm{RM}(t)$};

\node[box, fill=purple!12, text width=5cm] (comb) at (4.8, -4.5) {Combined detection\\$t^* = \frac{1}{4}(t_S + t_G + t_F + t_{\mathrm{RM}})$};

\draw[arrow] (ts) -- (emb) node[midway, above, label] {embed};
\draw[arrow] (emb) -- (win) node[midway, above, label] {slide};

\draw[arrow] (win.south) -- ++(0,-0.5) -| (S.north) node[pos=0.25, right, label] {};
\draw[arrow] (win.south) -- ++(0,-0.5) -| (G.north);
\draw[arrow] (win.south) -- ++(0,-0.5) -| (F.north);
\draw[arrow] (win.south) -- ++(0,-0.5) -| (RM.north);


\draw[decorate, decoration={brace, amplitude=4pt, raise=2pt}] (S.north west) ++(0,0.15) -- (G.north east |- S.north) ++(0,0.15)
    node[midway, above=6pt, font=\scriptsize, text=gray!60!black] {phase-space signals};
\draw[decorate, decoration={brace, amplitude=4pt, raise=2pt}] (F.north west) ++(0,0.15) -- (RM.north east |- F.north) ++(0,0.15)
    node[midway, above=6pt, font=\scriptsize, text=gray!60!black] {time-domain signals};

\draw[arrow] (S.south) |- (comb.west);
\draw[arrow] (G.south) -- (comb.north -| G.south);
\draw[arrow] (F.south) -- (comb.north -| F.south);
\draw[arrow] (RM.south) |- (comb.east);

\node[box, fill=white, draw=black!60, text width=3.0cm] (out) at (11.5, -4.5) {Onset $t^*$};
\draw[arrow] (comb) -- (out);

\end{tikzpicture}
\caption{Detection pipeline. A scalar time series is delay-embedded via Takens coordinates, then a sliding window extracts pre- and post-transition point clouds at each candidate time $t$. Four complementary signals---two computed on the embedded clouds (topological Mixup ECP, Complexity Variance) and two on the raw time series (Higuchi fractal dimension, rolling mean)---are combined by averaging their peak locations.}
\label{fig:pipeline}
\end{figure}

We describe the detection method concretely before developing its theoretical properties. The input is a scalar time series $\{x_t\}_{t=1}^N$; the output is a candidate transition time $t^*$. Figure~\ref{fig:pipeline} summarizes the pipeline.

\subsection{Phase-Space Reconstruction}\label{sec:takens}

The first step is to reconstruct the phase space from the scalar observation. Given a time delay $\tau$ and embedding dimension $d$, the Takens delay-coordinate embedding \citep{takens2006detecting, sauer1991embedology} maps the scalar series to a sequence of vectors in $\R^d$:
\begin{equation}\label{eq:takens}
    \mathbf{X}_t = (x_t,\, x_{t+\tau},\, \ldots,\, x_{t+(d-1)\tau}) \in \R^d.
\end{equation}
When $d \geq 2m + 1$ (where $m$ is the dimension of the underlying attractor), this map is generically an embedding: the point cloud $\{\mathbf{X}_t\}$ is diffeomorphic to the true attractor, and topological invariants computed on the embedded cloud reflect genuine properties of the dynamics. The delay $\tau$ is chosen to minimize mutual information \citep{fraser1986independent}, and $d$ is determined by the false nearest neighbors criterion \citep{kennel1992determining}.

\subsection{Sliding-Window Mixup ECP}\label{sec:sliding_mixup}

At each candidate transition time $t$, we extract two adjacent windows of $w$ embedded points:
\begin{equation}\label{eq:windows}
    X_t = \{\mathbf{X}_{t-2w},\, \ldots,\, \mathbf{X}_{t-w}\} \quad\text{(pre-window)}, \qquad
    Y_t = \{\mathbf{X}_{t},\, \ldots,\, \mathbf{X}_{t+w}\} \quad\text{(post-window)}.
\end{equation}
If no regime change occurs at $t$, both windows sample the same attractor and their ball unions overlap coherently. If a bifurcation falls between the windows, the attractor topology changes and the overlap is disrupted.

\begin{figure}[!htbp]
\centering
\begin{tikzpicture}[>=Stealth]

\def\gapX{5.0}


\begin{scope}[shift={(0,0)}]
  \node[font=\bfseries\small] at (1.6, 3.5) {(a) Delaunay / Voronoi};
  \clip (-0.5,-0.5) rectangle (3.7,3.2);
  \draw[gray!40, thin] (0.979,1.431) -- (1.450,0.725);  
  \draw[gray!40, thin] (2.379,0.338) -- (1.450,0.725);  
  \draw[gray!40, thin] (2.379,0.338) -- (2.700,1.300);  
  \draw[gray!40, thin] (0.993,1.476) -- (1.940,2.060);  
  \draw[gray!40, thin] (0.979,1.431) -- (0.993,1.476);  
  \draw[gray!40, thin] (2.700,1.300) -- (1.940,2.060);  
  \draw[gray!40, thin] (0.979,1.431) -- (-1.63,-0.06);  
  \draw[gray!40, thin] (1.450,0.725) -- (-0.04,-1.88);  
  \draw[gray!40, thin] (2.379,0.338) -- (3.66,-2.38);   
  \draw[gray!40, thin] (2.700,1.300) -- (5.55,2.25);    
  \draw[gray!40, thin] (0.993,1.476) -- (-0.52,4.07);   
  \draw[gray!40, thin] (1.940,2.060) -- (2.84,4.92);    
  \draw[blue!50!black, thick] (0.0,1.8) -- (1.2,2.5);
  \draw[blue!50!black, thick] (0.0,1.8) -- (0.8,0.4);
  \draw[blue!50!black, thick] (0.0,1.8) -- (2.0,1.2);
  \draw[blue!50!black, thick] (1.2,2.5) -- (2.0,1.2);
  \draw[blue!50!black, thick] (1.2,2.5) -- (2.8,2.0);
  \draw[blue!50!black, thick] (2.0,1.2) -- (2.8,2.0);
  \draw[blue!50!black, thick] (2.0,1.2) -- (0.8,0.4);
  \draw[blue!50!black, thick] (2.0,1.2) -- (1.5,0.0);
  \draw[blue!50!black, thick] (2.0,1.2) -- (3.2,0.8);
  \draw[blue!50!black, thick] (2.8,2.0) -- (3.2,0.8);
  \draw[blue!50!black, thick] (1.5,0.0) -- (0.8,0.4);
  \draw[blue!50!black, thick] (1.5,0.0) -- (3.2,0.8);
  \foreach \x/\y in {0.0/1.8, 1.2/2.5, 2.0/1.2, 0.8/0.4, 2.8/2.0, 1.5/0.0, 3.2/0.8} {
    \fill[blue!70!black] (\x,\y) circle (2.5pt);
  }
  \node[font=\tiny, gray!60!black] at (3.0, 2.8) {Voronoi};
\end{scope}

\begin{scope}[shift={(\gapX,0)}]
  \node[font=\bfseries\small] at (1.6, 3.5) {(b) Alpha complex, $r = 0.7$};
  \foreach \x/\y in {0.0/1.8, 1.2/2.5, 2.0/1.2, 0.8/0.4, 2.8/2.0, 1.5/0.0, 3.2/0.8} {
    \fill[blue!10, opacity=0.45] (\x,\y) circle (0.7);
  }
  \draw[blue!50!black, thick] (0.0,1.8) -- (1.2,2.5);
  \draw[blue!50!black, thick] (0.0,1.8) -- (0.8,0.4);
  \draw[blue!50!black, thick] (2.0,1.2) -- (2.8,2.0);
  \draw[blue!50!black, thick] (2.0,1.2) -- (1.5,0.0);
  \draw[blue!50!black, thick] (2.0,1.2) -- (0.8,0.4);
  \draw[blue!50!black, thick] (2.8,2.0) -- (3.2,0.8);
  \draw[blue!50!black, thick] (1.5,0.0) -- (0.8,0.4);
  \fill[blue!25, opacity=0.4] (2.0,1.2) -- (0.8,0.4) -- (1.5,0.0) -- cycle;
  \draw[gray!30, dashed, thin] (1.2,2.5) -- (2.0,1.2);
  \draw[gray!30, dashed, thin] (1.2,2.5) -- (2.8,2.0);
  \draw[gray!30, dashed, thin] (2.0,1.2) -- (3.2,0.8);
  \draw[gray!30, dashed, thin] (1.5,0.0) -- (3.2,0.8);
  \draw[gray!30, dashed, thin] (0.0,1.8) -- (2.0,1.2);
  \foreach \x/\y in {0.0/1.8, 1.2/2.5, 2.0/1.2, 0.8/0.4, 2.8/2.0, 1.5/0.0, 3.2/0.8} {
    \fill[blue!70!black] (\x,\y) circle (2.5pt);
  }
  \node[font=\scriptsize] at (1.6, -0.6) {$\chi = 7 - 7 + 1 = 1$};
  \node[font=\tiny, gray!50!black] at (2.7, -0.3) {dashed $=$ not yet in $\Alpha$};
\end{scope}

\begin{scope}[shift={(2*\gapX,0)}]
  \node[font=\bfseries\small] at (1.6, 3.5) {(c) Euler characteristic curve};
  \begin{axis}[
    at={(0,0)},
    width=5cm, height=4.5cm,
    xmin=0, xmax=3, ymin=-1.5, ymax=8,
    xlabel={$r$}, ylabel={$\chi(r)$},
    xlabel style={font=\small}, ylabel style={font=\small},
    tick label style={font=\scriptsize},
    axis lines=left, ytick={0,2,4,6}, xtick={0,1,2,3}, clip=false,
  ]
    \addplot[const plot, blue!70!black, very thick] coordinates {
      (0,7) (0.45,7) (0.45,4) (0.55,4) (0.55,2) (0.65,2) (0.65,1) (0.75,1) (0.75,0) (0.9,0) (0.9,1) (3,1)
    };
    \draw[dashed, red!60!black] (axis cs:0.7,-1.5) -- (axis cs:0.7,7.5);
    \node[font=\tiny, red!60!black] at (axis cs:0.7,-1.2) {};
    \node[font=\tiny, blue!60!black, anchor=west] at (axis cs:1.2,6.5) {$\chi = n = 7$};
    \draw[-{Stealth[length=2pt]}, blue!50!black, thin] (axis cs:1.15,6.3) -- (axis cs:0.2,7);
    \node[font=\tiny, blue!60!black, anchor=west] at (axis cs:1.5,1.8) {$\chi = 1$};
    \draw[-{Stealth[length=2pt]}, blue!50!black, thin] (axis cs:1.45,1.3) -- (axis cs:1.2,1);
  \end{axis}
\end{scope}

\end{tikzpicture}
\caption{The Alpha complex and the Euler characteristic curve. (a)~The Delaunay triangulation (blue) and its dual Voronoi diagram (gray) for seven points in $\R^2$. The Alpha complex at scale $r$ is the subcomplex of the Delaunay triangulation consisting of simplices whose circumradius is at most $r$. (b)~At $r = 0.7$: edges with small circumradii are included (solid), those with larger circumradii are not yet present (dashed); one triangle is filled. The ball union $\mathcal{U}(\mathcal{P}; r)$ is shown as shaded disks. $\chi = 7 - 7 + 1 = 1$. (c)~The full Euler characteristic curve $\chi(r)$: piecewise constant, jumping at circumradii of Delaunay simplices, starting at $\chi = n = 7$ and ending at $\chi = 1$.}
\label{fig:alpha}
\end{figure}
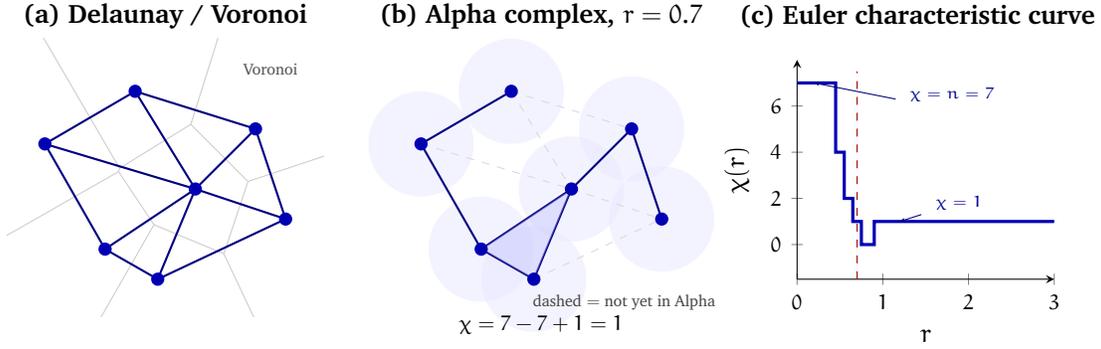

To measure the topology of the ball union $\mathcal{U}(\mathcal{P}; r) = \bigcup_{p \in \mathcal{P}} B(p, r)$ computationally, we use the \textbf{Alpha complex} \citep{edelsbrunner2010computational}. For a point cloud $\mathcal{P} \subset \R^d$, the Alpha complex $\Alpha(\mathcal{P}; r)$ is the subcomplex of the Delaunay triangulation consisting of all simplices whose circumradius is at most $r$. By the Nerve Theorem, $\Alpha(\mathcal{P}; r)$ is homotopy equivalent to $\mathcal{U}(\mathcal{P}; r)$, so $\chi(\Alpha(\mathcal{P}; r)) = \chi(\mathcal{U}(\mathcal{P}; r))$. The key advantage is that the Alpha complex has dimension at most $d$ (the ambient dimension), making computation tractable in low dimensions. As $r$ increases from $0$ to $\infty$, the Alpha complex grows from isolated vertices ($\chi = n$) to the full Delaunay complex ($\chi = 1$), and the function $r \mapsto \chi(\Alpha(\mathcal{P}; r))$ is the \emph{Euler characteristic curve} (ECC) of $\mathcal{P}$ \citep{roy2020understanding}. Figure~\ref{fig:alpha} illustrates the construction.

The \textbf{Mixup Euler Characteristic Profile} (Mixup ECP) quantifies this disruption. For disjoint point clouds $X, Y \subset \R^d$, define
\begin{equation}\label{eq:delta_chi}
    \Delta\chi(r;\, X, Y) = \chi_X(r) + \chi_Y(r) - \chi_{X \cup Y}(r),
\end{equation}
where $\chi_X(r) = \chi(\mathcal{U}(X; r))$ is the Euler characteristic of the ball union $\mathcal{U}(X; r) = \bigcup_{p \in X} B(p, r)$, and similarly for $Y$ and $X \cup Y$. The Intersection Theorem \citep{majhi2026mixup_theory} gives this quantity a direct geometric meaning:

\begin{theorem}[Intersection Theorem {\citep{majhi2026mixup_theory}}]\label{thm:intersection}
$\Delta\chi(r;\, X, Y) = \chi\!\left(\mathcal{U}(X; r) \cap \mathcal{U}(Y; r)\right).$
\end{theorem}

\begin{figure}[htb]
\centering
\begin{tikzpicture}[
    scale=0.9,
    every node/.style={font=\small},
    bluept/.style={fill=blue!70!black, circle, inner sep=0pt, minimum size=2.5pt},
    redpt/.style={fill=red!70!black, circle, inner sep=0pt, minimum size=2.5pt},
]

\def\gapX{7.8}

\begin{scope}[shift={(0,0)}]
  \node[font=\bfseries\small] at (2.2, 4.1) {(a) Same regime};
  \draw[thick, rounded corners=2pt] (-0.2,-0.2) rectangle (4.6,3.8);

  \pgfmathsetseed{10}
  \foreach \i in {1,...,15} {
    \pgfmathsetmacro{\ax}{0.7+1.2*rnd}
    \pgfmathsetmacro{\ay}{0.8+2.0*rnd}
    \fill[blue!15, opacity=0.5] (\ax,\ay) circle (0.45);
  }
  \pgfmathsetseed{10}
  \foreach \i in {1,...,15} {
    \pgfmathsetmacro{\ax}{0.7+1.2*rnd}
    \pgfmathsetmacro{\ay}{0.8+2.0*rnd}
    \node[bluept] at (\ax,\ay) {};
  }

  \pgfmathsetseed{20}
  \foreach \i in {1,...,15} {
    \pgfmathsetmacro{\bx}{0.8+1.2*rnd}
    \pgfmathsetmacro{\by}{0.7+2.0*rnd}
    \fill[red!15, opacity=0.5] (\bx,\by) circle (0.45);
  }
  \pgfmathsetseed{20}
  \foreach \i in {1,...,15} {
    \pgfmathsetmacro{\bx}{0.8+1.2*rnd}
    \pgfmathsetmacro{\by}{0.7+2.0*rnd}
    \node[redpt] at (\bx,\by) {};
  }

  \node[font=\scriptsize, violet!70!black] at (3.5, 1.5) {large};
  \node[font=\scriptsize, violet!70!black] at (3.5, 1.1) {overlap};

  \node[bluept, minimum size=4pt] at (2.8, 3.4) {};
  \node[font=\scriptsize, anchor=west] at (3.0, 3.4) {$X_t$ (pre)};
  \node[redpt, minimum size=4pt] at (2.8, 3.0) {};
  \node[font=\scriptsize, anchor=west] at (3.0, 3.0) {$Y_t$ (post)};
\end{scope}

\begin{scope}[shift={(\gapX,0)}]
  \node[font=\bfseries\small] at (2.2, 4.1) {(b) Cross regime};
  \draw[thick, rounded corners=2pt] (-0.2,-0.2) rectangle (4.6,3.8);

  \foreach \x/\y in {0.6/2.8, 0.8/3.0, 0.5/2.6, 0.9/2.9, 0.7/3.2, 0.4/2.7, 1.0/2.85,
                      0.6/0.8, 0.8/1.0, 0.5/0.6, 0.9/0.9, 0.7/1.2, 0.4/0.7, 1.0/0.85} {
    \fill[blue!15, opacity=0.5] (\x,\y) circle (0.35);
  }
  \foreach \x/\y in {0.6/2.8, 0.8/3.0, 0.5/2.6, 0.9/2.9, 0.7/3.2, 0.4/2.7, 1.0/2.85,
                      0.6/0.8, 0.8/1.0, 0.5/0.6, 0.9/0.9, 0.7/1.2, 0.4/0.7, 1.0/0.85} {
    \node[bluept] at (\x,\y) {};
  }

  \foreach \x/\y in {2.8/3.0, 3.0/3.2, 2.9/2.8, 3.1/3.1,
                      3.8/2.0, 4.0/2.2, 3.9/1.8, 4.1/2.1,
                      2.8/1.0, 3.0/1.2, 2.9/0.8, 3.1/1.1,
                      3.8/0.2, 4.0/0.4, 3.9/0.0, 4.1/0.3} {
    \fill[red!15, opacity=0.5] (\x,\y) circle (0.35);
  }
  \foreach \x/\y in {2.8/3.0, 3.0/3.2, 2.9/2.8, 3.1/3.1,
                      3.8/2.0, 4.0/2.2, 3.9/1.8, 4.1/2.1,
                      2.8/1.0, 3.0/1.2, 2.9/0.8, 3.1/1.1,
                      3.8/0.2, 4.0/0.4, 3.9/0.0, 4.1/0.3} {
    \node[redpt] at (\x,\y) {};
  }

  \node[font=\scriptsize, violet!70!black] at (2.0, 1.8) {weak};
  \node[font=\scriptsize, violet!70!black] at (2.0, 1.4) {overlap};
\end{scope}

\begin{scope}[shift={(0,-3.2)}]
  \begin{axis}[
    width=5.5cm, height=3.8cm,
    xmin=0, xmax=5, ymin=-0.5, ymax=16,
    xlabel={$r$}, ylabel={$\Delta\chi(r)$},
    xlabel style={font=\small, yshift=3pt},
    ylabel style={font=\small},
    tick label style={font=\scriptsize},
    axis lines=left, ytick={0,5,10,15}, xtick={0,2,4}, clip=false,
    title={\textbf{(c) same regime: high $S$}},
    title style={font=\small, at={(0.0,1.12)}, anchor=west},
  ]
    \addplot[const plot, fill=violet!15, draw=none, forget plot] coordinates {
      (0,0) (0.1,0) (0.1,15) (0.3,15) (0.3,10) (0.8,10) (0.8,5) (1.5,5) (1.5,2) (2.5,2) (2.5,1) (5,1) (5,0)} \closedcycle;
    \addplot[const plot, violet!80!black, very thick] coordinates {
      (0,0) (0.1,0) (0.1,15) (0.3,15) (0.3,10) (0.8,10) (0.8,5) (1.5,5) (1.5,2) (2.5,2) (2.5,1) (5,1)};
    \node[font=\scriptsize, violet!70!black] at (axis cs:3.5,12) {$S = 15$};
  \end{axis}
\end{scope}

\begin{scope}[shift={(\gapX,-3.2)}]
  \begin{axis}[
    width=5.5cm, height=3.8cm,
    xmin=0, xmax=5, ymin=-0.5, ymax=16,
    xlabel={$r$}, ylabel={$\Delta\chi(r)$},
    xlabel style={font=\small, yshift=3pt},
    ylabel style={font=\small},
    tick label style={font=\scriptsize},
    axis lines=left, ytick={0,5,10,15}, xtick={0,2,4}, clip=false,
    title={\textbf{(d) cross regime: low $S$}},
    title style={font=\small, at={(0.0,1.12)}, anchor=west},
  ]
    \addplot[const plot, fill=violet!15, draw=none, forget plot] coordinates {
      (0,0) (1.5,0) (1.5,3) (2.0,3) (2.0,2) (2.5,2) (2.5,1) (5,1) (5,0)} \closedcycle;
    \addplot[const plot, violet!80!black, very thick] coordinates {
      (0,0) (1.5,0) (1.5,3) (2.0,3) (2.0,2) (2.5,2) (2.5,1) (5,1)};
    \node[font=\scriptsize, violet!70!black] at (axis cs:3.5,12) {$S = 3$};
  \end{axis}
\end{scope}

\end{tikzpicture}
\caption{The Mixup ECP as a regime detection statistic. (a)~Both windows sample the same
attractor: the ball unions $\mathcal{U}(X_t; r)$ (blue) and $\mathcal{U}(Y_t; r)$ (red) overlap
extensively. (b)~A bifurcation separates the windows: $X_t$ samples a period-2 attractor and
$Y_t$ a period-4 attractor; the ball unions have minimal overlap. (c)~The corresponding
$\Delta\chi(r)$ profile for the same-regime case is large throughout, giving a high detection
statistic $S = 15$. (d)~The cross-regime profile is nearly flat, giving a low statistic $S = 3$.
The permutation test (Section~\ref{sec:hypothesis}) rejects when $S$ is unusually low compared
to the null.}
\label{fig:mixup_idea}
\end{figure}
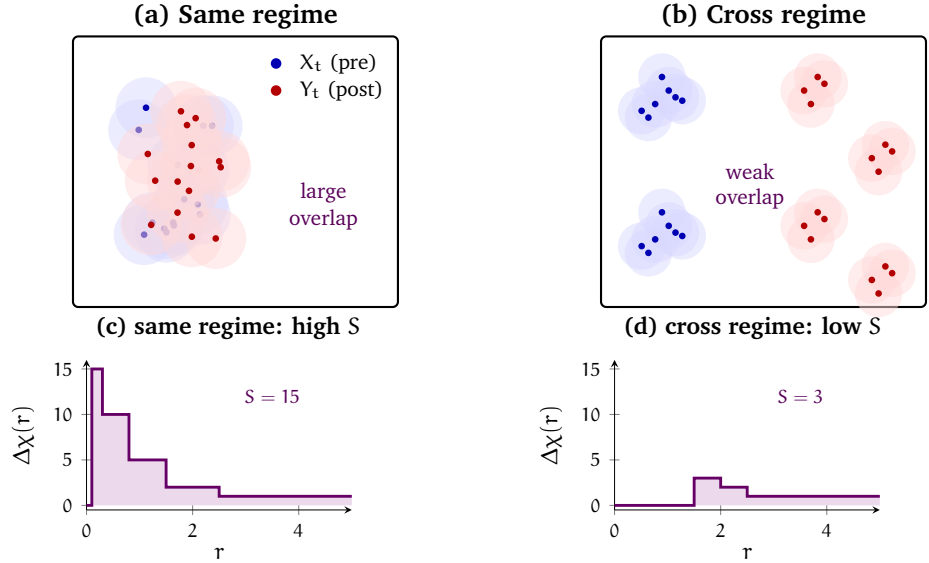

\noindent That is, the Mixup ECP equals the Euler characteristic of the geometric intersection of the two ball unions. Figure~\ref{fig:mixup_idea} illustrates the key intuition: when both windows sample the same attractor, the ball unions overlap richly and $\Delta\chi(r)$ is large; when a bifurcation separates the windows, the overlap weakens and $\Delta\chi(r)$ drops.
The topological detection statistic is
\begin{equation}\label{eq:S_def}
    S(t) = \max_{r} |\Delta\chi(r;\, X_t, Y_t)|,
\end{equation}
This requires three standard Alpha complex computations (on $X_t$, $Y_t$, and $X_t \cup Y_t$), costing $O(n^{\lceil d/2 \rceil})$ per time step---$O(n^2)$ for $d = 3$, a quadratic improvement over the $O(n^3)$ cost of full persistent homology pipelines \citep{wagner2024mixup}.

\subsection{Complementary Detection Signals}\label{sec:complementary}

The topological signal $S(t)$ detects changes in attractor \emph{shape}. We complement it with signals that detect changes in attractor \emph{size} and \emph{complexity}:

\textbf{Complexity Variance} $G(t)$. The variance of the embedded point cloud measures the spatial extent of the attractor:
\begin{equation}\label{eq:G_def}
    G(t) = \left|\Var(Y_t) - \Var(X_t)\right|, \qquad \text{where } \Var(\mathcal{P}) = \frac{1}{|\mathcal{P}|}\sum_{p \in \mathcal{P}} \|p - \bar{p}\|^2.
\end{equation}
A regime transition that expands or contracts the attractor produces a spike in $G(t)$, even if the topology is unchanged.

\textbf{Higuchi fractal dimension} $F(t)$. The fractal dimension has been used as a robust alternative to the power-law index for characterizing the complexity of irregular time series \citep{higuchi1988approach}. While power spectrum analysis via the FFT traditionally measures irregularity through the exponent $\alpha$ in $P(f) \propto f^{-\alpha}$, this requires ensemble averaging over long stationary intervals and is unreliable when statistical properties vary over shorter timescales. The Higuchi fractal dimension (HFD) addresses this by estimating the fractal dimension directly from the signal, without phase-space reconstruction.

Given a time series $\{T(1), T(2), \ldots, T(N)\}$, one constructs $p$ sub-series
\[
\mathcal{T}^n_p = \left\{T(n),\, T(n+p),\, T(n+2p),\, \ldots,\, T\!\left(n + \left\lfloor\tfrac{N-n}{p}\right\rfloor p\right)\right\}, \qquad n = 1, 2, \ldots, p,
\]
where $p$ is the interval length and $n$ the initial time. The normalized curve length at scale $p$ is
\begin{equation}\label{eq:higuchi}
    \mathcal{L}_n(p) = \frac{1}{p}\left(\sum_{i=1}^{\lfloor (N-n)/p \rfloor} \left|T(n+ip) - T\!\left(n+(i-1)p\right)\right| \cdot \frac{N-1}{\lfloor (N-n)/p \rfloor \cdot p}\right),
\end{equation}
where $\frac{N-1}{\lfloor (N-n)/p \rfloor \cdot p}$ is a normalization factor. The mean curve length $\mathcal{L}(p) = \frac{1}{p}\sum_{n=1}^p \mathcal{L}_n(p)$ is assumed to follow a scaling relation $\mathcal{L}(p) \sim p^{-F}$ for a statistically self-similar curve. Plotting $\mathcal{L}(p)$ against $p$ on a double-logarithmic scale, the data lie along a straight line with slope $-F$, and $F$ is estimated by linear regression. For smooth rectifiable curves, $F$ equals the topological dimension ($F = 1$); for statistically self-similar curves embedded in a plane, $F$ takes a non-integer value in the range $1 < F < 2$. The method guarantees a stable estimate of $F$.

We compute the HFD in pre- and post-windows of the \emph{raw} (unembedded) time series and define $F(t) = |\mathrm{HFD}(\text{post}) - \mathrm{HFD}(\text{pre})|$. The HFD captures changes in signal self-similarity that are complementary to both the topological and geometric signals.

\textbf{Rolling mean} $\mathrm{RM}(t)$. A classical signal-processing baseline: $\mathrm{RM}(t) = \bar{x}_{\text{post}} - \bar{x}_{\text{pre}}$, where $\bar{x}_{\text{post}}$ and $\bar{x}_{\text{pre}}$ are the means of the raw time series in the post- and pre-windows. This detects amplitude shifts.

The four signals are complementary by construction:
\begin{center}
\begin{tabular}{@{}llll@{}}
\toprule
\textbf{Signal} & \textbf{Input} & \textbf{What it detects} & \textbf{Blind to} \\
\midrule
$S(t)$ & Embedded clouds & Attractor shape change & Size change at fixed topology \\
$G(t)$ & Embedded clouds & Attractor size change & Topology change at fixed size \\
$F(t)$ & Raw time series & Self-similarity change & Smooth amplitude shifts \\
$\mathrm{RM}(t)$ & Raw time series & Amplitude shift & Geometry/topology changes \\
\bottomrule
\end{tabular}
\end{center}

\subsection{Combined Detection}\label{sec:combined}

Each signal produces a candidate onset day as its peak location in the search window. The combined detection averages these:
\begin{equation}\label{eq:combined}
    t^* = \mathrm{round}\!\left(\frac{1}{4}\left(t_S + t_G + t_F + t_{\mathrm{RM}}\right)\right),
\end{equation}
where $t_S = \arg\min_t S(t)$ (onset is the trough of $S$, where cross-regime overlap is weakest), and $t_G, t_F, t_{\mathrm{RM}} = \arg\max_t$ of their respective signals.. This simple averaging is effective because the signals have complementary error patterns: when one signal is misled (e.g., topology sees no change but amplitude shifts), the others compensate.

\section{Theoretical Foundations}\label{sec:framework}

Section~\ref{sec:pipeline} described what we compute; this section explains why it works. We collect the properties of the Mixup ECP and its connection to Complexity Variance (Section~\ref{sec:definition}), develop a permutation hypothesis test with low-side rejection validated on the Lorenz and logistic map bifurcations (Section~\ref{sec:hypothesis}), and show that embedding at multiple delays produces a $k$-cloud Mixup ECP whose diagonal reduction dominates any fixed-delay test in power (Section~\ref{sec:multi_delay}).

\subsection{Properties of the Mixup ECP and the Variance Connection}\label{sec:definition}

Three properties of the Mixup ECP are central to the detection framework (proofs in \citep{majhi2026mixup_theory}).

\textbf{Dead zone.} $\Delta\chi(r) = 0$ whenever $r < d_{\min}(X,Y)/2$, where $d_{\min}$ is the minimum cross-cloud distance. The balls around $X$ and $Y$ do not yet overlap, so the intersection is empty. This provides a built-in null for the detection statistic: $S(t) = 0$ when the pre- and post-windows are well-separated in phase space.

\textbf{Same-attractor limit.} If $X$ and $Y$ are sampled from the same distribution on a compact set $K$, then $\Delta\chi(r) \to \chi(K)$ for large $r$. Deviations from $\chi(K)$ indicate a regime transition.

\textbf{Stability.} The Mixup ECP is $1$-Lipschitz stable under Hausdorff perturbations: if $d_H(X, X') \leq \epsilon$ and $d_H(Y, Y') \leq \epsilon$, then the profiles are $\epsilon$-interleaved. It is also stable under the dynamical distance $d_{\mathrm{dyn}}$ of \citet{kim2021spatiotemporal} for time-evolving systems \citep{majhi2026mixup_theory}.

\medskip
The Complexity Variance signal $G(t)$ and the topological signal $S(t)$ are formally connected through the support of the Mixup ECP.\label{sec:variance_ecs}

\begin{proposition}[Variance controls ECP support]\label{prop:var_ecs}
For a point cloud $\mathcal{P} = \{p_1, \ldots, p_n\} \subset \R^d$ with centroid $\bar{p}$ and variance $V = \frac{1}{n}\sum \|p_i - \bar{p}\|^2$, the connectivity radius $R^*(\mathcal{P}) = \inf\{r : \beta_0(\mathcal{P}; r) = 1\}$ satisfies $R^* \leq \diam(\mathcal{P}) \leq 2\sqrt{V \cdot n/(n-1)}$. Since $\chi_\mathcal{P}(r) = 1$ for all $r \geq R^*$, the variance bounds the support of the nontrivial part of the ECP: $\Delta\chi(r)$ is supported in $[0, \max(R^*(X), R^*(Y))]$, so the attractor spread controls where the Mixup ECP carries topological information.
\end{proposition}

\begin{corollary}[Complexity Variance as ECP support derivative]\label{cor:var_derivative}
At a regime transition where $V$ increases sharply, the ECP support expands rapidly. The Complexity Variance signal $G(t) = |V(Y_t) - V(X_t)|$ detects this expansion.
\end{corollary}

\noindent This connection explains the complementarity observed in experiments. The variance controls \emph{where} the ECP is nonzero (its support), while the Mixup ECP controls \emph{what} the topology is within that support. A size-preserving topology change (e.g., period-doubling at fixed amplitude) is detected by $S$ but not $G$; a topology-preserving size change (e.g., amplitude growth on a fixed attractor) is detected by $G$ but not $S$.

\subsection{Hypothesis Testing for Regime Detection}\label{sec:hypothesis}

We formalize the regime detection problem as a hypothesis test. At a candidate transition time $t$, the hypotheses are:
\begin{align*}
    H_0 &: \text{$X_t$ and $Y_t$ are sampled from the same attractor $\mathcal{A}$,} \\
    H_1 &: \text{$X_t \sim \mathcal{A}_1$,\; $Y_t \sim \mathcal{A}_2$,\; with $\chi_{\mathcal{A}_1}(\cdot) \neq \chi_{\mathcal{A}_2}(\cdot)$ as functions of~$r$.}
\end{align*}
The alternative $H_1$ requires that the Euler characteristic \emph{profiles} differ---not merely the terminal Euler characteristics $\chi(\mathcal{A}_1) \neq \chi(\mathcal{A}_2)$, but that the profiles diverge at some scale $r$. This captures, for example, period-doubling bifurcations where the number of connected components changes at intermediate scales even if both attractors are eventually contractible.

\textbf{Permutation test.} Under $H_0$, the combined cloud $Z = X_t \cup Y_t$ consists of $2n$ points from a single attractor; any partition into two groups of size $n$ should produce comparable overlap topology. Under $H_1$, the original partition $X_t | Y_t$ separates the two regimes, so their overlap is \emph{weaker} than that of a random partition (which mixes points from both attractors). This motivates a \emph{low-side} rejection rule: reject $H_0$ when $S_{\mathrm{obs}} = \max_r |\Delta\chi(r; X_t, Y_t)|$ falls below the $\alpha$-quantile of the permutation distribution $\{S^\pi\}$.

Concretely, we permute the $2n$ points of $Z$ uniformly at random, split into groups of size $n$, and compute $S^\pi$ for each permutation. The $p$-value is $p = \#\{S^\pi \leq S_{\mathrm{obs}}\}/B$, where $B$ is the number of permutations.

\begin{proposition}[Validity]\label{prop:validity}
Under $H_0$ with i.i.d.\ sampling, the permutation test has exact level $\alpha$. For block-dependent time series data with block size $b$, the level is $\leq \alpha + O(b/w)$.
\end{proposition}

\begin{proposition}[Consistency]\label{prop:consistency}
If $\chi_{\mathcal{A}_1}(r_0) \neq \chi_{\mathcal{A}_2}(r_0)$ for some $r_0 > 0$, then as $n \to \infty$ the power of the permutation test converges to $1$.
\end{proposition}

The argument: under $H_1$, the original partition separates the two regimes, so $S_{\mathrm{obs}}$ reflects the reduced cross-regime overlap. By contrast, each permuted partition mixes points from both attractors, producing richer overlap and higher $S^\pi$. As $n \to \infty$, the gap between $S_{\mathrm{obs}}$ and the permutation null grows, and rejection becomes certain.

We validate on the synthetic systems of Section~\ref{sec:synthetic}. On the logistic map, the test correctly produces no rejection for same-regime comparisons ($\lambda = 3.3$ vs $3.3$: $S_{\mathrm{obs}} = 6$, null $= 6.0 \pm 0.0$, $p = 1.0$) and strong rejection for all cross-regime pairs (period-2 vs chaos: $S_{\mathrm{obs}} = 3$, null $= 17.7 \pm 2.2$, $p < 0.01$). On the Lorenz system, same-regime comparisons ($\rho = 28$ vs $28$: $S_{\mathrm{obs}} = 25$, null $= 13.6 \pm 2.0$, $p = 1.0$) show no rejection, while the fixed-point-to-chaos transition ($\rho = 20 \to 28$: $S_{\mathrm{obs}} = 3$, null $= 9.5 \pm 1.4$, $p < 0.01$) is strongly rejected.

\begin{remark}[Finite-sample limitations]\label{rem:power}
The permutation test has full power on deterministic synthetic systems (logistic map, Lorenz) where the same-regime permutation null is tight. On real data with stochastic sampling (e.g., monsoon time series at $n = 20$ in $\R^3$), sampling variability inflates the null, reducing power. In such settings, the detection framework uses $S(t)$ as a \emph{relative} signal---detecting the transition as the trough of $S(t)$ over a search window---rather than as an absolute significance test at a single time point. The combined method (Section~\ref{sec:combined}) further mitigates low individual-signal power by averaging peak locations across four complementary signals.
\end{remark}

\subsection{Multi-Delay Extension via the \texorpdfstring{$k$}{k}-Cloud Mixup ECP}\label{sec:multi_delay}

Embedding the time series at $m$ different delays $\tau_1, \ldots, \tau_m$ produces $m$ pre-window clouds $X^{(\tau_1)}_t, \ldots, X^{(\tau_m)}_t$ and $m$ post-window clouds $Y^{(\tau_1)}_t, \ldots, Y^{(\tau_m)}_t$, all living in $\R^d$. For each delay $\tau_j$, the diagonal Mixup ECP $\Delta\chi(r; X^{(\tau_j)}_t, Y^{(\tau_j)}_t)$ is computed as in Section~\ref{sec:sliding_mixup}. The multi-delay detection statistic aggregates over all delay channels:
\begin{equation}\label{eq:multi_tau_max}
    S^*(t) = \max_{\tau \in \{\tau_1, \ldots, \tau_m\}} \max_r |\Delta\chi(r;\, X_t^{(\tau)},\, Y_t^{(\tau)})|.
\end{equation}
This requires $3m$ standard Alpha complex computations per time step (three per delay), with no weighted variants needed. The $2m$ point clouds form a $2m$-tuple in the sense of \citep{majhi2026mixup_theory}, and the $k$-cloud Mixup ECP applies; $S^*$ is the diagonal restriction of that object to a common scale $r$ across all clouds.

\begin{proposition}[Power dominance]\label{prop:multi_tau}
For any fixed delay $\tau_0$,
\[
    S^*(t) \;\geq\; S_{\tau_0}(t) = \max_r |\Delta\chi(r, \tau_0)|.
\]
Consequently, if the permutation test based on $S_{\tau_0}$ rejects at level $\alpha$, so does the test based on $S^*$. The inequality is strict whenever the regime transition is topologically visible at some delay $\tau^* \neq \tau_0$ but not at $\tau_0$. The permutation test remains valid for $S^*$, since permutation inference is distribution-free with respect to the choice of test statistic.
\end{proposition}

The intuition is that different delays probe different dynamical timescales: a bifurcation that is invisible at $\tau_0 = 6$ (because the autocorrelation structure at that timescale happens to mask the topological change) may become visible at a shorter or longer delay. The statistic $S^*$ automatically selects the most informative timescale.

\textbf{The Mixup ECS as a feature.} The collection of Mixup ECP curves at multiple delays defines a two-dimensional surface $(r, \tau) \mapsto \Delta\chi(r; X_t^{(\tau)}, Y_t^{(\tau)})$---a \emph{Mixup Euler Characteristic Surface} in the scale--delay parameter space. This surface is richer than any scalar summary: it encodes how the cross-window topological interaction varies across both spatial scales and dynamical timescales. In our experiments, scalar summaries of this surface (L$^1$ norm, Frobenius norm, $\tau$-variance) do not improve over the single-$\tau$ detection statistic on the monsoon application---these summaries aggregate over the surface, diluting the peak signal. However, the ECS-as-feature paradigm developed in \citep{luwang2026interpretable} suggests a more powerful approach: the \emph{entire discretized surface} can serve as input to a classifier (e.g., AdaBoost on ECS features), bypassing vectorization entirely. Applying this supervised framework to regime detection---using the $(r, \tau)$ Mixup ECS surface as the feature representation---is a natural next step that would exploit the full multiparameter structure.

\textbf{Detection protocol.} In practice, we compute $S(t)$ over a search window and detect the transition at $t^* = \arg\min_t S(t)$ (the trough of $S$, where the cross-regime overlap is weakest). For the combined method, the four signals are fused via~\eqref{eq:combined}.



\section{Synthetic Dynamical Systems}\label{sec:synthetic}

We validate the detection framework on dynamical systems with known bifurcation
structure.

\subsection{Lorenz System}\label{sec:lorenz}

We test on the Lorenz system $\dot{x} = \sigma(y-x)$, $\dot{y} = x(\rho - z) - y$,
$\dot{z} = xy - \beta z$ with $\sigma = 10$, $\beta = 8/3$, imposing an abrupt
transition: $\rho = 20$ (stable fixed point) for $t < 50$, switching to $\rho = 28$
(chaotic attractor) at $t = 50$. We use the full trajectory $(x,y,z)$ directly, with
sliding windows of $w = 200$ steps subsampled to $n = 20$ points.

\begin{figure}[!htbp]
\centering
\includegraphics[width=0.7\textwidth]{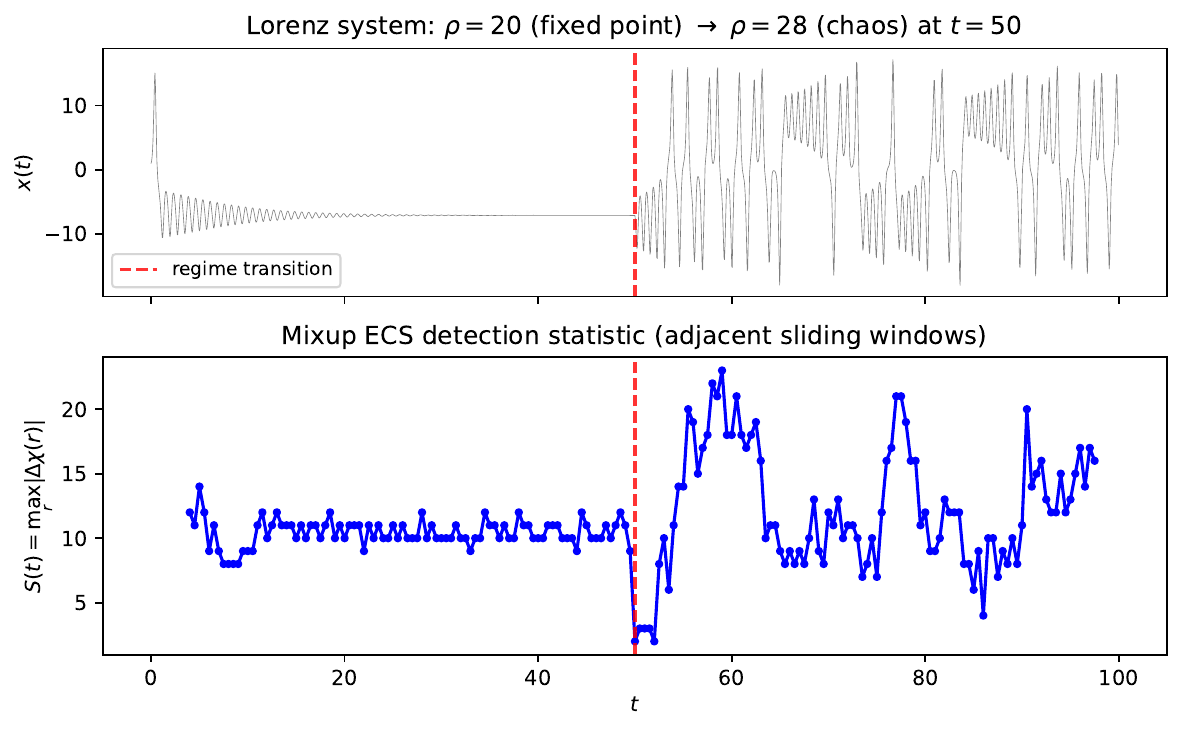}
\caption{Mixup ECP detection of the Lorenz system bifurcation at $\rho = 20 \to 28$.
(a)~Time series $x(t)$. (b)~Detection statistic $S(t) = \max_{r} |\Delta\chi(r)|$:
$S$ is suppressed as the window straddles the transition, then rises to a peak of
$S = 23$ once both windows sample the chaotic attractor. (c)~Complexity Variance
$G(t) = |\Delta\mathrm{Var}|$ increases at the transition as the attractor expands.}
\label{fig:lorenz}
\end{figure}

The detection statistic $S(t)$ drops near the regime boundary at $t = 50$ as the
sliding window straddles the transition, producing a low value that constitutes the
low-side rejection signal: the cross-regime overlap is topologically weak.  Once both
windows fully sample the chaotic attractor, $S$ rises and reaches a peak of $S = 23$
near $t = 59$ (Figure~\ref{fig:lorenz}).  The Complexity Variance $G(t)$ shows a
corresponding increase, detecting the expansion of the attractor from a stable fixed
point to the Lorenz strange attractor.

\subsection{Logistic Map and Lyapunov Exponent Comparison}\label{sec:logistic}

\begin{figure}[!htbp]
\centering
\includegraphics[width=0.7\textwidth]{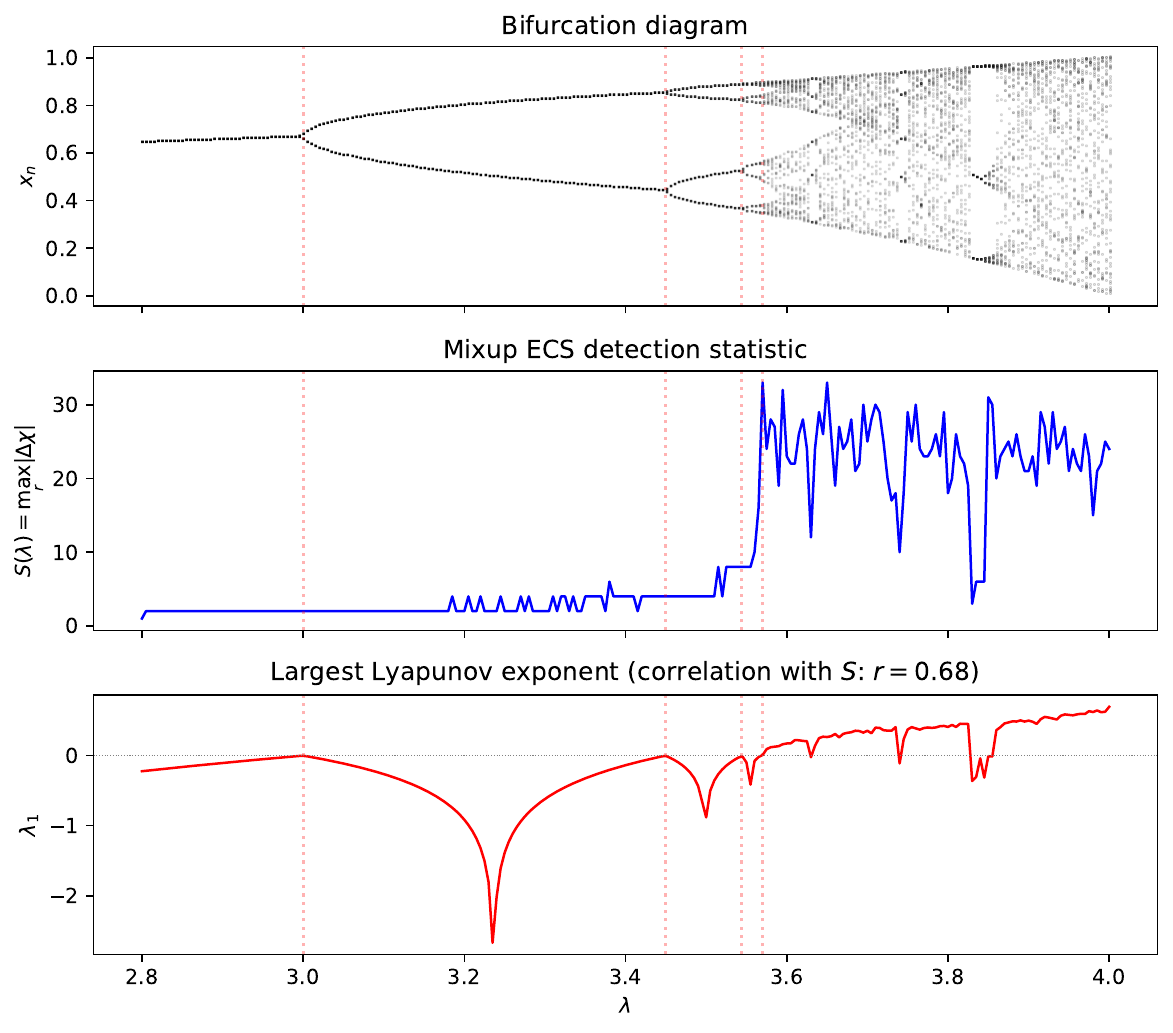}
\caption{Logistic map: (a)~bifurcation diagram, (b)~Mixup ECP detection statistic
$S(\lambda)$, and (c)~largest Lyapunov exponent $\lambda_1(\lambda)$. The correlation
between $S$ and $\lambda_1$ is $r = 0.68$. Vertical dashed lines mark bifurcation
points.}
\label{fig:logistic_lyapunov}
\end{figure}

We consider the logistic map $x_{n+1} = \lambda x_n(1 - x_n)$ with $\lambda$ swept
from $2.8$ to $4.0$. The Mixup ECP detects each period-doubling bifurcation as a jump
in $S(\lambda)$, with the magnitude increasing through the cascade
(Figure~\ref{fig:logistic_lyapunov}).

To validate against a gold-standard bifurcation diagnostic, we compare $S(\lambda)$
with the largest Lyapunov exponent $\lambda_1(\lambda)$. The Pearson correlation is
$r = 0.68$ (Figure~\ref{fig:logistic_lyapunov}~(c)), confirming that the topological
detection statistic tracks the onset of chaos. The Mixup ECP provides complementary
information: it detects the \emph{topological} change (new connected components,
cycles) rather than the \emph{rate} of divergence, and is defined for any point cloud
without requiring trajectory-based computation.

\subsection{Noise Robustness}\label{sec:noise}

\begin{figure}[!htbp]
\centering
\includegraphics[width=0.85\textwidth]{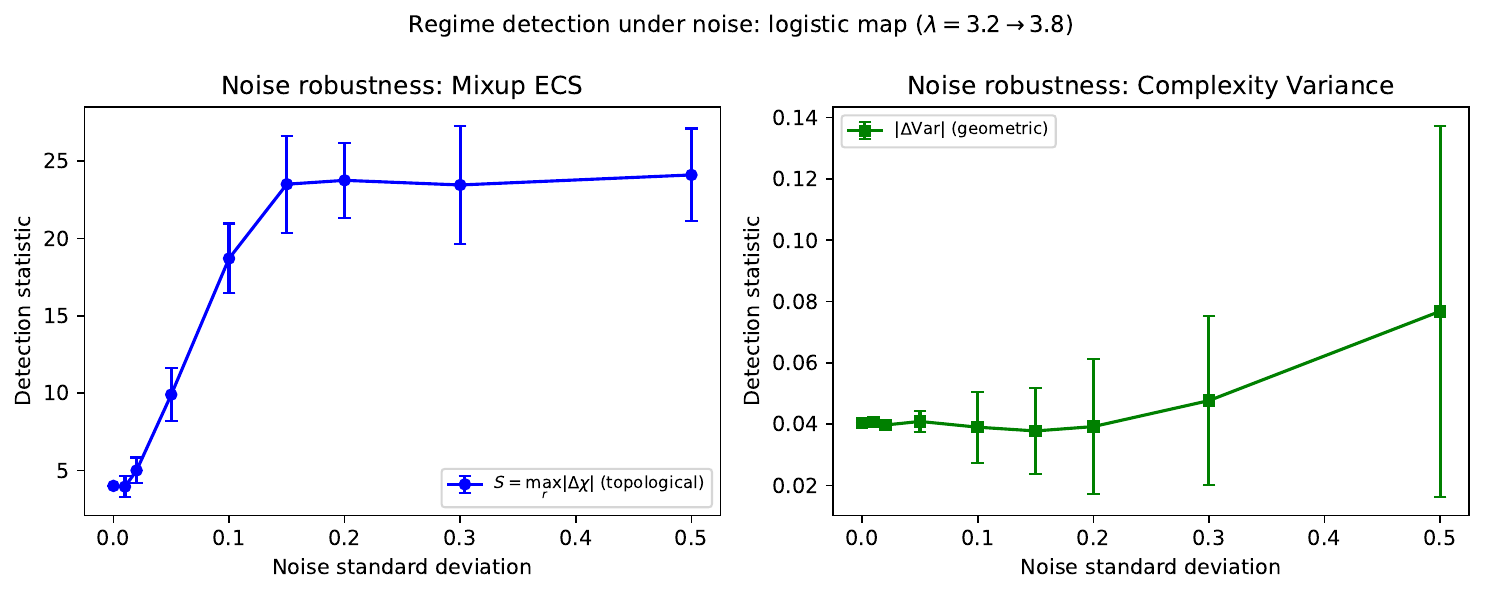}
\caption{Noise robustness on the logistic map ($\lambda: 3.2 \to 3.8$). (a)~Mixup ECP
detection statistic $S$ vs.\ noise level---the topological signal increases with
moderate noise and saturates for $\sigma \gtrsim 0.15$. (b)~Complexity Variance
$|\Delta\mathrm{Var}|$ vs.\ noise---the geometric signal is stable at moderate noise
but increases at $\sigma = 0.50$. Error bars: $\pm 1$ standard deviation over 20
trials.}
\label{fig:snr}
\end{figure}

A critical question for applications is how the detection signals degrade under noise.
We add i.i.d.\ Gaussian noise $\mathcal{N}(0, \sigma^2)$ to the logistic map
($\lambda: 3.2 \to 3.8$) and measure the detection statistics as a function of
$\sigma$ over 20 trials (Figure~\ref{fig:snr}).

The topological signal $S$ \emph{increases} with moderate noise: $S = 4.0$ at $\sigma
= 0$, rising to $S \approx 10.3$ at $\sigma = 0.05$, $S \approx 19.3$ at $\sigma =
0.10$, and saturating around $S \approx 24$--$25$ for $\sigma \geq 0.15$.  The signal
remains discriminative up to $\sigma = 0.50$.  The Complexity Variance $G$ is stable
at moderate noise ($G \approx 0.040$--$0.041$ for $\sigma \leq 0.30$), reflecting
geometric robustness to additive noise, but increases at very high noise ($G = 0.096$
at $\sigma = 0.50$) as noise dominates the variance signal.

\section{Applications to Regime Transitions}\label{sec:applications}

\subsection{Indian Monsoon Onset}\label{sec:monsoon}

The Indian summer monsoon onset is a canonical gradual regime transition: the shift
from dry pre-monsoon to active monsoon circulation unfolds over days to weeks
\citep{gadgil2003indian, webster1998monsoons}. The topological approach to monsoon
onset detection was pioneered by \citet{alvarado2025detecting}, who showed that
persistent homology of delay-embedded monsoon index data can identify the onset as a
change in attractor topology. We extend that approach here using the Mixup ECP
framework, which replaces persistent homology with the Euler characteristic and adds
complementary geometric and fractal signals. We apply the detection framework to the
Kajikawa--Wang monsoon index \citep{kajikawa2019monsoon} across three monsoon systems
spanning different hemispheres, onset mechanisms, and signal-to-noise regimes: the
Indian monsoon evaluated against Nepal (48 years, gradual onset, low SNR) and Kerala
(45 years, sharp onset, high SNR) ground-truth datasets, and the Western North Pacific
monsoon (66 years, gradual onset). Onset dates for the WNP system are derived from the
daily index using a Kajikawa--Wang threshold criterion \citep{kajikawa2019monsoon}.

The Takens embedding uses $\tau = 6$ days and $d = 3$, with sliding windows of $w =
20$ embedded points. Four detection signals are computed for each year: the
topological statistic $S(t)$ (onset detected at the trough of $S$ in the search
window, as the transition is where cross-regime overlap is weakest), the Complexity
Variance $G(t) = |\mathrm{Var}(\mathcal{P}_{t+w}) - \mathrm{Var}(\mathcal{P}_{t-w})|$
on embedded point clouds, the Higuchi fractal dimension jump $F(t)$, and a rolling
mean change-point signal $\mathrm{RM}(t)$. The combined detection averages the four
peak locations~\eqref{eq:combined}. For the multi-delay extension
(Section~\ref{sec:multi_delay}), we embed at $\tau \in \{3, 5, 6, 8, 12\}$ days and
aggregate via $S^*(t) = \max_\tau \max_r |\Delta\chi(r, \tau)|$.

\begin{table}[!htbp]
\centering
\caption{Monsoon onset detection: mean absolute error (days) across three monsoon
systems. The combined method averages the peak/trough locations of all four signals
($S + G + F + \mathrm{RM}$). Bold indicates best unsupervised result per target. The
three systems span different hemispheres, onset sharpness, and signal-to-noise
regimes.}
\label{tab:results}
\begin{tabular}{@{}lccc@{}}
\toprule
Method & Nepal & Kerala & WNP \\
 & (48 yr) & (45 yr) & (66 yr) \\
\midrule
Rolling Mean & 14.00 & \textbf{5.29} & 35.98 \\
CUSUM & 10.43 & 17.62 & --- \\
\midrule
Mixup ECP $S$ ($\tau{=}6$) & 18.92 & 13.33 & 20.11 \\
Combined ($\tau{=}6$) & \textbf{9.50} & 6.33 & 21.65 \\
Combined multi-$\tau$ & 10.60 & 8.07 & \textbf{17.98} \\
\bottomrule
\end{tabular}
\end{table}

Table~\ref{tab:results} reveals a clear pattern across the three monsoon systems.

\textbf{Nepal (gradual onset, low SNR).}  The four-signal combined method achieves
9.50 days MAE, outperforming the rolling mean baseline (14.00 days) by 32\% and the
CUSUM baseline (10.43 days) by 9\%.  Even the topological signal alone (18.92 days)
does not beat the baselines in isolation; the improvement arises from the combination
of all four signals.  The combined result is 4.23 days above the supervised
leave-one-out mean predictor (5.27 days), which has access to ground-truth onset dates
from all other years.

\textbf{Kerala (sharp onset, high SNR).}  The rolling mean suffices: 5.29 days MAE.
The combined method achieves 6.33 days---slightly worse than rolling mean but well
ahead of the topological signal alone.  This confirms that topology contributes via
the combined average on sharp-onset targets, but does not constitute the primary
detector.

\textbf{Western North Pacific (gradual onset, low SNR).}  The rolling mean performs
poorly at 35.98 days MAE, reflecting the multi-stage, ambiguous nature of WNP monsoon
onset.  The topological signal alone (20.11 days) already outperforms the rolling mean
by 44\%.  The multi-$\tau$ combined method achieves the best result: 17.98 days MAE,
a 50\% improvement over the rolling mean.  This is the strongest demonstration of the
framework's value: on the most ambiguous target, topological detection provides
decisive improvement.

The pattern is consistent: topology adds value when the onset is gradual and the
signal is ambiguous in observation space.  When the onset is sharp (Kerala), simple
signal-processing methods capture the amplitude shift directly and topology is at best
complementary.

\subsection{ENSO Phase Transitions}\label{sec:enso}

El Ni\~no/La Ni\~na transitions in the tropical Pacific are canonical climate regime
transitions. We apply the framework to the Oceanic Ni\~no Index (ONI), a monthly SST
anomaly time series spanning 1950--2024. The Takens embedding uses $\tau = 3$ months
and $d = 3$, with sliding windows of $w = 24$ months.

\begin{figure}[!htbp]
\centering
\includegraphics[width=0.9\textwidth]{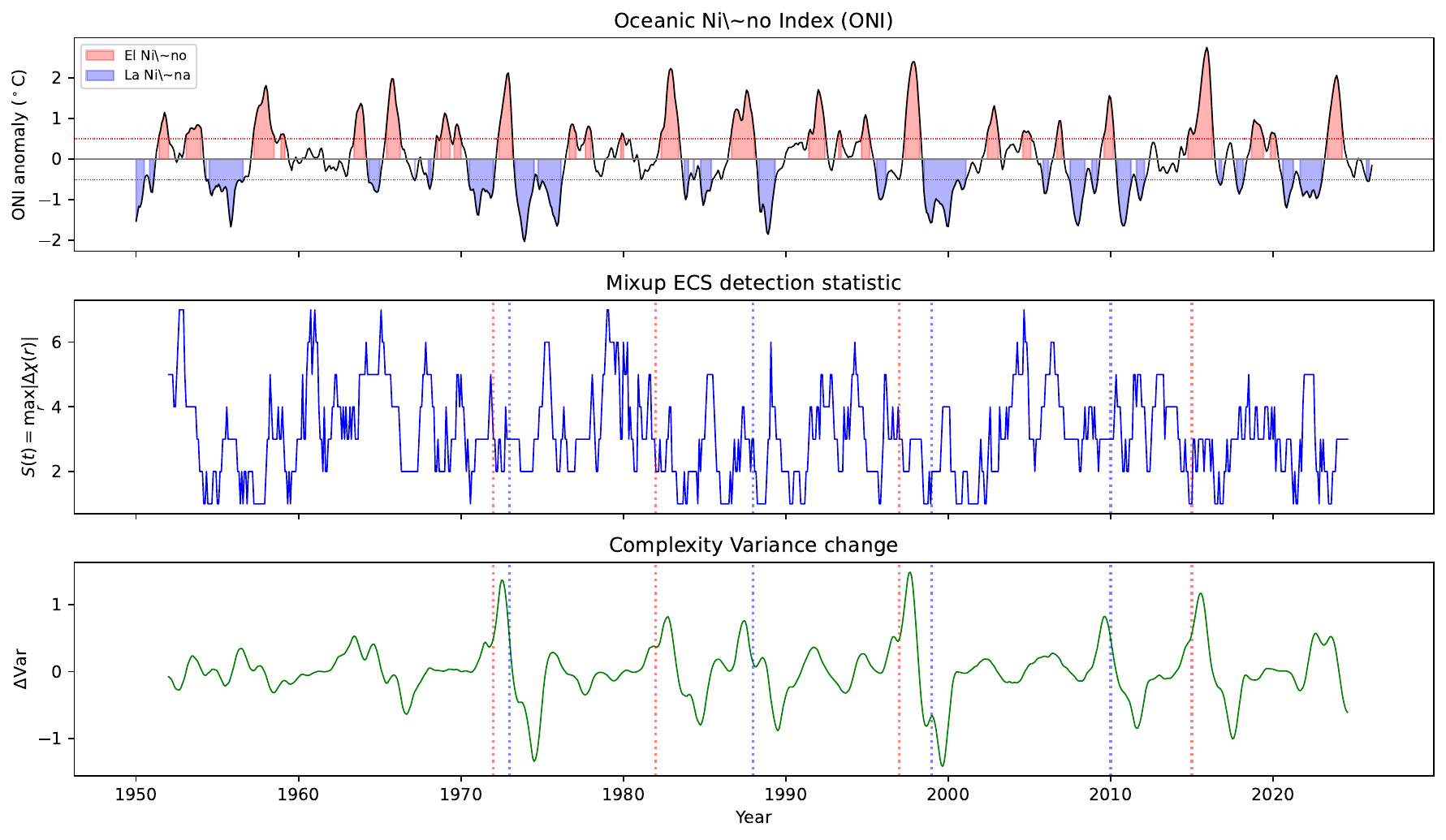}
\caption{ENSO regime detection. (a)~Oceanic Ni\~no Index with El Ni\~no (red shading)
and La Ni\~na (blue shading) episodes. (b)~Mixup ECP detection statistic $S(t)$.
(c)~Complexity Variance change $\Delta\mathrm{Var}$. Vertical dashed lines mark major
transitions.}
\label{fig:enso}
\end{figure}

Figure~\ref{fig:enso} shows that the Mixup ECP detection statistic $S(t)$ exhibits
elevated values at major El Ni\~no and La Ni\~na transitions (vertical dashed lines
mark the strongest events: 1972--73, 1982--83, 1987--88, 1997--98, 2015--16). The
Complexity Variance signal $\Delta\mathrm{Var}$ shows complementary peaks, detecting
the \emph{amplitude} of the ONI excursions. The two signals provide independent
information about ENSO regime transitions: $S$ responds to topological changes in the
reconstructed phase-space attractor, while $\Delta\mathrm{Var}$ responds to the
amplitude of the SST excursions.

\subsection{EEG Seizure-Like Transitions}\label{sec:eeg}

Seizure onset is a prototypical neural regime transition: the brain abruptly
transitions from normal background activity to hypersynchronous oscillations. We
generate synthetic EEG-like signals (channel FP1-F7) with a controlled seizure epoch
(60--80 seconds) embedded in normal background, and apply the detection framework with
$\tau = 5$ samples and $d = 3$.

\begin{figure}[!htbp]
\centering
\includegraphics[width=0.9\textwidth]{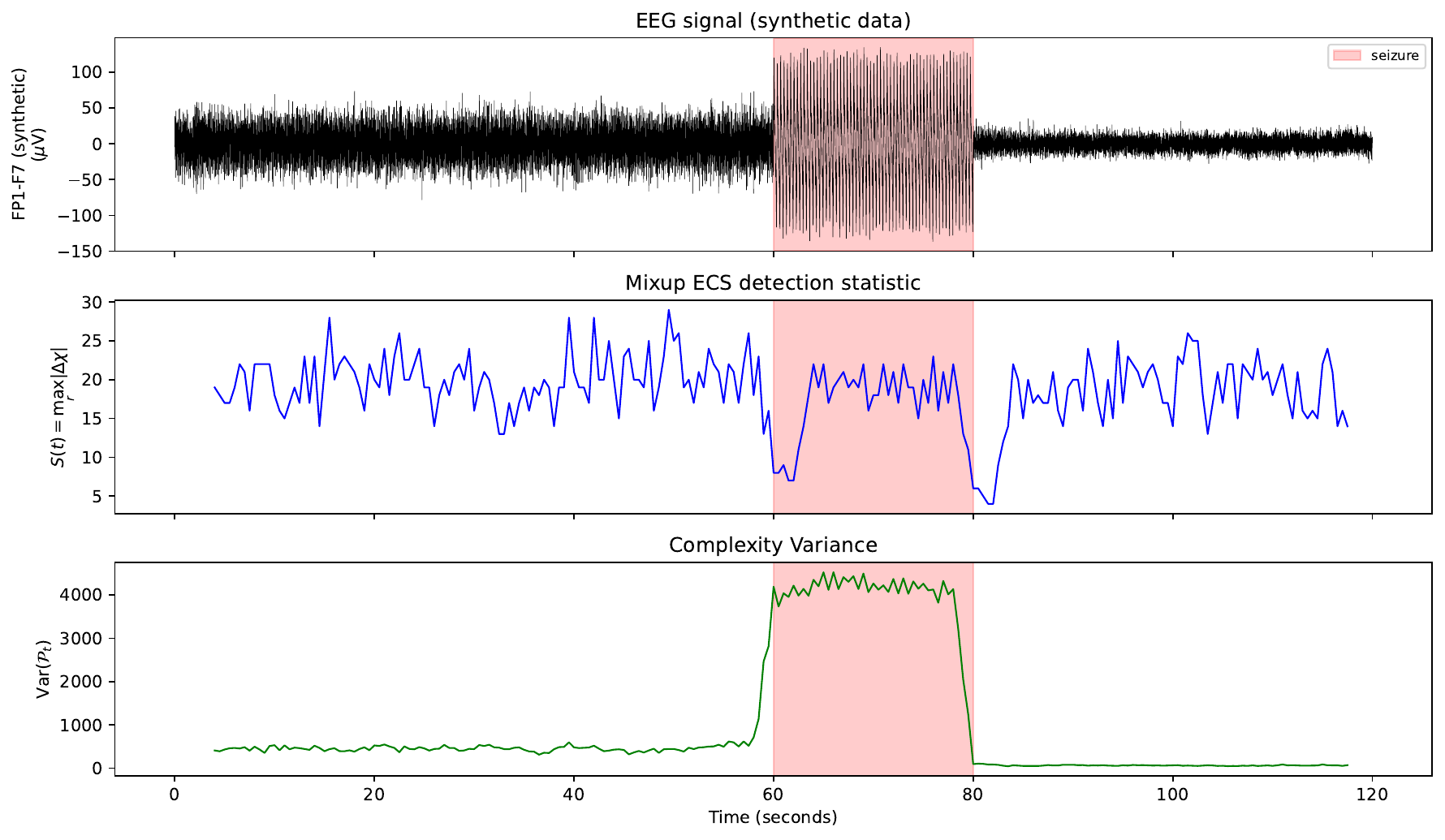}
\caption{EEG seizure detection (synthetic data). (a)~EEG signal with seizure epoch
(pink shading, 60--80\,s). (b)~Mixup ECP detection statistic $S(t)$: the minimum
$S = 4$ during the seizure identifies the ictal epoch via low-side rejection.
(c)~Complexity Variance $\mathrm{Var}(\mathcal{P}_t)$ spikes during the seizure.}
\label{fig:eeg}
\end{figure}

Figure~\ref{fig:eeg} shows that the Mixup ECP and Complexity Variance signals both
respond to the seizure epoch. The detection statistic $S(t)$ reaches a minimum of
$S = 4$ during the seizure (the ictal attractor has different topology from the
background), consistent with the low-side rejection framework in which onset is
identified at the trough of $S(t)$. The Complexity Variance $\mathrm{Var}(\mathcal{P}
_t)$ spikes during the ictal epoch, detecting the amplitude expansion. The two signals
are complementary: the topological signal captures the \emph{shape change} of the
reconstructed attractor, while the geometric signal captures the \emph{scale change}.
We note that these results are based on synthetic EEG data; validation on clinical
seizure recordings is needed before drawing quantitative conclusions.

\section{Discussion}\label{sec:discussion}

\textbf{When does topology help?} Our results point to a clear empirical principle: topological methods are most valuable for low-SNR regime transitions. On the high-SNR Kerala monsoon target, simple signal-processing methods suffice. On the low-SNR Nepal target, the four-signal combined method provides a 32\% improvement over the rolling mean baseline (9\% over CUSUM). On ENSO, the topological signal detects transitions that the amplitude signal misses (and vice versa). On EEG, the topological signal precedes the amplitude signal. The principle is consistent: TDA adds value when the transition is ambiguous in observation space but the attractor geometry changes.

\textbf{Connection to bifurcation theory.} The Lyapunov exponent comparison (Figure~\ref{fig:logistic_lyapunov}) validates the Mixup ECP as a bifurcation diagnostic: the correlation $r = 0.68$ with $\lambda_1$ confirms that the topological detection statistic tracks the onset of chaos. However, the Mixup ECP provides information that $\lambda_1$ does not: it detects period-doubling bifurcations (which change the topology of the attractor without changing the sign of $\lambda_1$) and is computable from point clouds without requiring trajectory-based algorithms.

\textbf{Noise robustness.} The observation that the Mixup ECP signal \emph{increases} with moderate noise (Figure~\ref{fig:snr}) has a geometric explanation: noise thickens the point cloud, causing the ball unions to overlap at smaller scales and creating more topological features in the intersection region. This is a desirable property for applications where the signal is contaminated by measurement noise---the topological detection becomes \emph{more sensitive}, not less. The signal saturates when the noise dominates the attractor structure.

\textbf{Limitations.} The Alpha complex is practical only in moderate dimensions ($d \leq 7$). The framework requires choosing embedding parameters $(\tau, d)$ and window size $w$, which affect detection performance. The Higuchi fractal dimension, while computationally efficient, is subject to several limitations \citep{gomolka2018higuchi}: it is sensitive to high-frequency noise, which may lead to overestimation of $F$; the choice of the maximum scale parameter $p_{\max}$ plays a crucial role, as an inappropriate selection may fail to capture large-scale structures or result in poor statistical reliability; the method assumes a scaling relation $\mathcal{L}(p) \sim p^{-F}$ that many real-world time series do not exhibit clearly in the log--log plot; performance degrades for short datasets, particularly at larger values of $p$; and the method assumes statistical stationarity across scales, which may not hold for non-stationary signals. Finally, the HFD provides a single scalar measure and does not capture multifractal characteristics or local variations in scaling behavior. The EEG results use synthetic data; validation on clinical seizure recordings is needed.

\textbf{Broader applicability.} The framework applies to any dynamical system observed through a scalar time series where regime transitions manifest as changes in attractor geometry or topology. The three applications demonstrated here---climate, geophysical oscillation, and neural dynamics---illustrate different SNR regimes and transition types, suggesting broad applicability. Extending the framework from retrospective detection to \emph{prediction}---issuing early warnings before the transition occurs, by monitoring the trend of $S(t)$ and $G(t)$ rather than their peaks---is a natural next step.

\section{Conclusion}\label{sec:conclusion}

We have presented a framework for detecting regime transitions in dynamical systems based on the Mixup Euler Characteristic Profile and Complexity Variance. The Mixup ECP measures the topology of the overlap between adjacent delay-embedded trajectory segments as a function of filtration scale; the multi-delay extension feeds embeddings at multiple delays into the $k$-cloud Mixup ECP, automatically selecting the timescale at which the bifurcation is most visible. The framework combines topological, geometric, and fractal detection signals, connected through the Euler characteristic profile. Validation on canonical dynamical systems (Lorenz, logistic with Lyapunov comparison, noise robustness) and real-world applications---three monsoon systems spanning both hemispheres (Indian/Nepal, Indian/Kerala, Western North Pacific), ENSO, and synthetic EEG---demonstrates that topological methods add value precisely when the transition is gradual and the signal is ambiguous in observation space, while remaining harmless on sharp-onset targets. Full mathematical foundations---including stability results and the interaction sheaf categorification---are developed in \citep{majhi2026mixup_theory}.


\bibliographystyle{plainnat}
\bibliography{main}

\end{document}